\documentclass[%
 aip,
 amsmath,amssymb,
 reprint,%
]{revtex4-1}

\usepackage{graphicx}
\usepackage{dcolumn}
\usepackage{bm}

\usepackage{amssymb}
\usepackage{amsmath}
\usepackage{amsthm}
\usepackage[utf8]{inputenc}
\usepackage[T1]{fontenc}
\usepackage{mathptmx}
\usepackage{etoolbox}
\usepackage{xcolor}

\makeatletter
\def\@email#1#2{%
 \endgroup
 \patchcmd{\titleblock@produce}
  {\frontmatter@RRAPformat}
  {\frontmatter@RRAPformat{\produce@RRAP{*#1\href{mailto:#2}{#2}}}\frontmatter@RRAPformat}
  {}{}
}%
\makeatother
\begin{document}

\preprint{AIP/123-QED}

\title[On solutions to the continuum version of the Kuramoto model with identical oscillators]{On solutions to the continuum version of the Kuramoto model with identical oscillators}
\author{A. R. Krueger}
 \affiliation{Department of Mathematics, Virginia Tech, Blacksburg, Virginia}
\author{S. Rengaswami}%
\affiliation{ 
Department of Mathematics, University of Tennessee, Knoxville, Tennessee
}%

\author{ R. N. Leander}
\email{rachel.leander@mtsu.edu}
\affiliation{Department of Mathematical Sciences, Middle Tennessee State University, Murfreesboro, Tennessee}%

\date{\today}

\begin{abstract}
The Kuramoto model provides a concrete mathematical realization of emergent synchrony in a population of phase-coupled oscillators.
Since Kuramoto's publication, \textit{Oscillations, Waves, and Turbulence}, researchers have worked to better characterize solution dynamics. In this paper, we combine the method of characteristics with an iterative technique to prove existence and uniqueness of solutions to the continuum version of the Kuramoto model in the special case where oscillators are identical and characterize the oscillator density in terms of a limiting projected characteristic. The model characterization, in turn, provides for global asymptotic analysis of the system via sub and super solutions. We believe the unique approach to model analysis developed here has the potential to yield novel results on other more complex versions of the extensively studied Kuramoto model.
\end{abstract}

\maketitle

\section{\label{sec:level1}Background}

The Kuramoto model was proposed in 1975 by Dr. Yoshiki Kuramoto as a concrete mathematical representation of emergent synchrony in a network of coupled oscillators \cite{KuraLecture}. More specifically, the Kuramoto model describes the synchronization of globally coupled phase oscillators, where coupling between oscillator pairs is proportional to the sine of their phase difference. Oscillator coupling is uniform in the sense that the constant of proportionality is uniform throughout the population. In addition, the population is characterized by a distribution of natural velocities. In \cite{kura} Kuramoto posited the existence of uniformly distributed (desynchronized) and travelling wave (partially-synchronized) asymptotic solutions for the continuum version of the model and discerned conditions for the synchronization-desynchronization transition in terms of the model coupling strength. He considered the nature of the bifurcation at the critical value for the coupling strength, however, the analysis proved difficult \cite{kura}. Additional results were obtained by Dr. John Crawford \cite{Letters}. Through the use of center manifold theory and linear analysis, Crawford studied the stability of solutions to the Kuramoto model \cite{crawford1}$^,$ \cite{crawford2}. Then, in 2008, Dr. Edward Ott and Dr. Thomas Antonsen developed an ansatz for model solutions in case the initial conditions satisfy a geometric Fourier series. They also developed an exact, closed form solution to the Kuramoto model given a Lorentzian distribution of natural velocities. Thus, Ott and Antonsen provided existence of solutions to the model under some assumptions on the initial data. Their ansatz is well-studied and has been used to develop numerical results in physical mathematics\cite{OttAntonsen}.

Here we study solutions to the continuum version of the Kuramoto model with identical oscillators. Our method differs from that of Ott and Antonsen in that it is based on the method of characteristics for solving first order partial differential equations. We characterize the oscillator density in terms of limiting characteristics, and characterize the asymptotic behavior of solutions, showing the synchronous state is globally, asymptotically stable. Although our results are specific to the case of identical natural velocities, we anticipate our methods will extend to the case of distributed natural velocities as well. \\

After writing this paper, we learned of closely related work by Chiba et al. \cite{Chiba1, Chiba2}. Indeed, the solution characterization presented here was previously derived by Chiba et al. as a weak solution formula \cite{Chiba1} in the more general case of distributed natural velocities (i.e., nonidentical oscillators). This solution formula was subsequently used to prove the Kuramoto conjecture \cite{Chiba2}, which includes a local stability analysis of the model's trivial and synchronous solutions.

\section{Model exposition}
In this paper we consider the continuum version of the Kuramoto model: 
{\fontsize{8}{9}\begin{align}
\frac{\partial\rho}{\partial t}(t,\theta,\omega)=&-\frac{\partial}{\partial \theta}\rho(t,\theta,\omega)\omega\\
&- \frac{\partial}{\partial\theta}\left(k \rho(t,\theta,\omega)\int_0^{2\pi} \int_{\mathbb{R}}{\rho(t,\hat{\theta},\hat{\omega}) \gamma(\hat{\omega})\sin(\hat{\theta}-\theta)d\hat{\omega}d\hat{\theta}}\right)\nonumber\\
\rho(0,\theta,\omega)=&\rho_0(\theta,\omega),
\end{align}}
where $(t,\theta,\omega)\in[0,\infty)\times\mathbb{R}\times\mathbb{R}$ and $\rho(t,\theta; \omega)$ is the $2\pi$-periodic (in $\theta$) probability density function for a population of oscillators, $\omega$ represents a natural velocity, $\gamma(\omega)$ is the probability density function of the oscillator natural velocities, $\theta$ represents an angle, $t$ represents time, and $k$ determines the coupling strength. In particular, 
\begin{equation}
\int_{a}^{b} \gamma(\omega)\mbox{ }d\omega
\end{equation}gives the probability that an oscillator has natural velocity between $a$ and $b$, and since $gamma$ is a probability density function, \begin{eqnarray}
\int_{-\infty}^{\infty} \gamma(\omega)\mbox{ }d\omega=1.\end{eqnarray}
Also,
\begin{equation}
\int_{-\infty}^{\infty} \rho(t,\theta,\omega)\gamma(\omega)\mbox{ }d\omega
\end{equation}gives the probability density of oscillator angles at time $t,$ so that for $0\leq a\leq b$, 
\begin{equation}
\int_{a}^{b}\int_{-\infty}^{\infty} \rho(t,\theta,\omega)\gamma(\omega)\mbox{ }d\omega\mbox{ }d\theta
\end{equation}gives the probability an oscillator has an angle between $a$ and $b$ at time $t$, and \begin{eqnarray}
\int_{0}^{2\pi}\int_{-\infty}^{\infty} \rho(t,\theta,\omega)\gamma(\omega)\mbox{ }d\omega\mbox{ }d\theta=1.\label{prob_density1}\end{eqnarray}

Note in addition, since the solution is 2$\pi$-periodic, the coupling term 
\begin{equation}\int_0^{2\pi} \int_{-\infty}^{\infty}{\rho(t,\hat{\theta},\omega) \gamma(\omega)\sin(\hat{\theta}-\theta)\mbox{ }d\omega\mbox{ }d\hat{\theta}}\label{coupling_term}
\end{equation}
is equal to
 \begin{equation}
 \int_{j2\pi}^{(j+1)2\pi} \int_{-\infty}^{\infty}{\rho(t,\hat{\theta},\omega) \gamma(\omega)\sin(\hat{\theta}-\theta)\mbox{ }d\omega\mbox{ }d\hat{\theta}}
\end{equation}
where $j2\pi\leq \theta\leq (j+1)2\pi,  j\in \mathbb{Z}$. Hence the model behaves as if oscillators are only coupled to those near them. 

In this first paper, we show existence and uniqueness of continuously differentiable solutions of a simplified model where all oscillators have a single natural velocity, i.e. $ \omega\equiv0 $. Under this assumption, the model conditions become
{\fontsize{8}{7}
\begin{align}
	\label{simplifiedmodel}\frac{\partial \rho}{\partial t} (t,\theta)&=-\frac{\partial}{\partial \theta}\left(k \rho(t,
\theta)\int_0^{2\pi}{\rho(t,\hat{\theta})\sin(\hat{\theta}-\theta)\mbox{ }d\hat{\theta}}\right),\\
\label{simplemodelIC}\rho(0,\theta)&=\rho_0(\theta),\\
\label{ProbabilityDensityWOmegaZero}\int_{0}^{2\pi} \rho_0(\theta)\mbox{ }d\theta&=1\\
\label{periodicity}\rho_0(\theta)&=\rho_0(\theta+2\pi),
\end{align}
}
where $(t,\theta)\in[0,\infty)\times\mathbb{R}$. Note that condition (\ref{ProbabilityDensityWOmegaZero}) stipulates that $\rho_0$ is a probability density function, and condition (\ref{periodicity}) stipulates that $\rho_0$ is $2\pi$-periodic. We would like $\rho$ to share these properties. 

It will be useful to write (\ref{simplifiedmodel}) as 
{\small
\begin{eqnarray} \label{PDE}
\frac{\partial \rho}{\partial t}(t,\theta)=&& -\frac{\partial \rho}{\partial \theta}(t,\theta)\left(k \int_0^{2\pi}{\rho(t,\hat{\theta})\sin(\hat{\theta}-\theta)d\hat{\theta}}\right)\nonumber\\
&&+\rho(t,\theta)\left(k \int_0^{2\pi}{\rho(t,\hat{\theta})\cos(\hat{\theta}-\theta)d\hat{\theta}}\right),
\end{eqnarray}
}
where we have differentiated through the integral in (\ref{simplifiedmodel}) with respect to $\theta$.

\section{Existence of Solutions}
We seek to show the existence of a nonnegative, global, solutions of (\ref{simplifiedmodel})-(\ref{periodicity}). In so doing, we assume the initial data, $\rho_0(\theta)$, is a nonnegative, $C^{2}$, $2\pi$-periodic probability density function. We employ an iterative method involving an approximating sequence, characterize the elements of the sequence, and establish convergence to a solution.  

We consider the following approximating sequence 
{\fontsize{8}{7}
\begin{eqnarray}\label{appxseq}
\frac{\partial \rho_{n+1}}{\partial t}(t,\theta)=&&-\frac{\partial \rho_{n+1}}{\partial \theta}(t,\theta)\left(k \int_0^{2\pi}{\rho_{n}(t,\hat{\theta})\sin(\hat{\theta}-\theta)d\hat{\theta}}\right)\nonumber\\
&&+\rho_{n+1}(t,\theta)\left(k \int_0^{2\pi}{\rho_{n}(t,\hat{\theta})\cos(\hat{\theta}-\theta)d\hat{\theta}}\right)
\end{eqnarray}
}
where $\rho_{n+1}$ satisfies the initial condition 
\begin{equation}
\label{IC_z}\rho_{n+1}(0,\theta)=\rho_0(\theta).
\end{equation}  First we show that if $\rho_n$ is $C^2$, non-negative, and $2\pi$-periodic, there exists a global solution to equation (\ref{appxseq}).  Note that (\ref{appxseq}) is a first order, linear, $2\pi$-periodic differential equation, hence we will employ the method of characteristics.  

Let
{
\begin{align}
F_n(x,y,z,p,q):=&p+q\left(k\int_0^{2\pi}{\rho_{n}(x,\hat{\theta})\sin(\hat{\theta}-y)d\hat{\theta}}\right)\\
&-z\left(k\int_0^{2\pi}{\rho_{n}(x,\hat{\theta})\cos(\hat{\theta}-y)d\hat{\theta}}\right),
\end{align}
}
so (\ref{appxseq}) can be expressed as

\begin{equation}
F_n\left(t,\theta ,\rho_{n+1},\frac{\partial \rho_{n+1}}{\partial t},\frac{\partial \rho_{n+1}}{\partial \theta}\right)=0.
\end{equation}


The characteristic equations associated with (\ref{appxseq}) are
{\fontsize{8}{7}
\begin{eqnarray}
	\label{xyzchar1}\dot{x}_{n+1}(s)&=&1,\\
	\label{xyzchar2}\dot{y}_{n+1}(s)&=&k\int_{0}^{2 \pi}\rho_{n}(x_{n+1}(s), \hat{\theta})\sin(\hat{\theta}-y_{n+1}(s))d\hat{\theta},\\
	\label{xyzchar3}\dot{z}_{n+1}(s)&=&z_{n+1}(s)k\int_{0}^{2 \pi}\rho_{n}(x_{n+1}(s), \hat{\theta})\cos(\hat{\theta}-y_{n+1}(s))d\hat{\theta},
\end{eqnarray}
}
together with
{\fontsize{8}{7}
\begin{eqnarray}
\label{p1p2char1}\dot{p}_{n+1}(s)=&& -q_{n+1}(s)k\int_0^{2\pi}{\frac{\partial}{\partial x}\rho_{n}(x_{n+1}(s),\hat{\theta})\sin(\hat{\theta}-y_{n+1}(s))d\hat{\theta}}\nonumber,\\
&&+ z_{n+1}(s)k\int_0^{2\pi}{\frac{\partial}{\partial x}\rho_{n}(x_{n+1}(s),\hat{\theta})\cos(\hat{\theta}-y_{n+1}(s))d\hat{\theta}}\\
&&+ p_{n+1}(s)k\int_0^{2\pi}{\rho_{n}(x_{n+1}(s),\hat{\theta})\cos(\hat{\theta}-y_{n+1}(s))d\hat{\theta}}\nonumber,\\
\label{p1p2char2}\dot{q}_{n+1}(s)=&& 2q_{n+1}(s)k\int_0^{2\pi}{\rho_{n}(x_{n+1}(s),\hat{\theta})\cos(\hat{\theta}-y_{n+1}(s))d\hat{\theta}}\nonumber\\
&&+z_{n+1}(s)k\int_0^{2\pi}{\rho_{n}(x_{n+1}(s),\hat{\theta})\sin(\hat{\theta}-y_{n+1}(s))d\hat{\theta}}.
\end{eqnarray}
}
In deriving the characteristic equations, we have used the following expressions for the partial derivatives of $F$ with respect to $x, y$ and $z$, which are found by differentiating through the associated integrals with respect to $x$ (using the fact that $\rho_n$ is $C^2$).
{\small
\begin{eqnarray}
(F_n)_x=&&q_{n+1}k\int_0^{2\pi}{\frac{\partial}{\partial x} \rho_{n}(x_{n+1},\hat{\theta})\sin(\hat{\theta}-y_{n+1})d\hat{\theta}}\nonumber\\
&&-z_{n+1}k\int_0^{2\pi}{\frac{\partial}{\partial x} \rho_{n}(x_{n+1},\hat{\theta})\sin(\hat{\theta}-y_{n+1})d\hat{\theta}},\\
(F_n)_y=&&-q_{n+1}k\int_0^{2\pi}{\rho_{n}(x_{n+1},\hat{\theta})\cos(\hat{\theta}-y_{n+1})d\hat{\theta}}\nonumber\\
&&-z_{n+1}k\int_0^{2\pi}{\rho_{n}(x_{n+1},\hat{\theta})\sin(\hat{\theta}-y_{n+1})d\hat{\theta}},\\
(F_n)_z=&&-k\int_0^{2\pi}{\rho_{n}(x_{n+1},\hat{\theta})\cos(\hat{\theta}-y_{n+1})d\hat{\theta}}.
\end{eqnarray}
}

We show that the system of equations (\ref{xyzchar1})-(\ref{xyzchar3}) has a unique global solution under any initial data.  However, for the purpose of solving (\ref{appxseq}) subject to (\ref{IC_z}), the  following initial data are of interest:
\begin{equation}
    (x_{n+1}(0), y_{n+1}(0), z_{n+1}(0))=(0,\theta_0,\rho_0(\theta_0))\label{specinitcondchar1-3}
\end{equation} 
Eventually we will need initial data for the $ p $ and $ q $ components as well, but these initial conditions will be specified later. For convenience, we will denote solutions of (\ref{xyzchar1})-(\ref{xyzchar3}) subject to initial data 
\begin{equation}(x_{n+1}(s_0), y_{n+1}(s_0), z_{n+1}(s_0))=(x_0,y_0,z_0)\label{geninitcondchar1-3}\end{equation} 
as 
{\fontsize{8}{7}
\begin{equation}\label{geninitcondchar1-3_solution}
X_{n+1}(s; s_0, X_0)=(x_{n+1}(s; s_0, X_0), y_{n+1}(s; s_0, X_0), z_{n+1}(s; s_0, X_0)),
\end{equation}}
where 
\begin{equation}
X_0=(x_0,y_0,z_0),
\end{equation}
and solutions of (\ref{xyzchar1})-(\ref{xyzchar3}) subject to (\ref{specinitcondchar1-3}) by   \begin{equation}
X_{n+1}(s;\theta_0)=(x_{n+1}(s;\theta_0),y_{n+1}(s;\theta_0),z_{n+1}(s;\theta_0)).
\end{equation}
\newtheorem*{lem1}{Lemma 1} 
\begin{lem1}
If $\rho_0(\theta)$ and $\rho_n(t,\theta)$ are $C^2$ in all arguments, then there exists a unique, global solution, $X_{n+1}(s;s_0,X_0)=(x_{n+1}(s;s_0,X_0),y_{n+1}(s;s_0,X_0),z_{n+1}(s;s_0,X_0)),$ to (\ref{xyzchar1})-(\ref{xyzchar3}) subject to the initial conditions of (\ref{geninitcondchar1-3}). Moreover, $X_{n+1}(s;s_0,X_0)$
is $C^2$ with respect to all arguments.  
\end{lem1}

\begin{proof}
Let $G_{n+1}(s,x,y,z)$ be defined by the right hand side of (\ref{xyzchar1})-(\ref{xyzchar3}). We will refer to the components of $ G_{n+1} $ as $ [G_{n+1}]_i $, $i=1, 2, 3$, and we will restrict the values of $ s $ to a compact interval $ [0,T] $.  Note that the characteristics equations can be solved in sequence. The unique, global solution to (\ref{xyzchar1}) satisfying the first component of (\ref{geninitcondchar1-3}) is 
\begin{equation} 
x_{n+1}=s+x_0,
\end{equation}
which is clearly $C^2$ in all arguments. Substituting $x_{n+1}(s)=s+x_0$ into (\ref{xyzchar2}) gives a differential equation for $y_{n+1}$ in $s$ and $y_{n+1}$ alone. Since sine is smooth and its derivative is bounded, $[G_{n+1}]_2$ is continuously differentiable with respect $y$ and its derivative with respect to $y$ is uniformly bounded for $(s,y)\in[0, T] \times \mathbb{R}$.
In particular, $[G_{n+1}]_2$ is uniformly Lipschitz continuous with respect to $y$ on $[0, T] \times \mathbb{R}.$ Hence, there exists a unique, $C^1$, global (i.e. defined for all $s\in \mathbb{R}$)  solution of $ (\ref{xyzchar1})- (\ref{xyzchar2}) $ satisfying the first two components of (\ref{geninitcondchar1-3}) (See Theorem 2.2, p.38 and Corollary 2.6, p.41 of \cite{teschl}).\\
With $x_{n+1}$ and $y_{n+1}$ in hand, $ (\ref{xyzchar3}) $ can be solved by separation of variables. This gives a continuously differentiable and global solution:
{\fontsize{9}{8}
\begin{eqnarray}
 \label{zsolgen}z_{n+1}(s; s_0, X_0)= z_0 e^{k\int_{s_0}^{s}\int_{0}^{2\pi}\rho_{n}(\tau,\hat{\theta})cos(\hat{\theta}-y_{n+1}(\tau;s_0,X_0))d\hat{\theta}\mbox{ }d\tau} .
 \end{eqnarray}
 }
 Note that $z_{n+1}(s, s_0, X_0)\neq{0}$ provided $z_0\neq 0$ and $z_{n+1}(s, s_0, X_0)\equiv 0$, provided $z_0=0$.  

Thus, unique $C^1$ functions $x_{n+1}(s;s_0,X_0), y_{n+1}(s;s_0,X_0)$ and $z_{n+1}(s;s_0,X_0)$ satisfying $(\ref{xyzchar1})-(\ref{xyzchar3})$, subject to $(\ref{geninitcondchar1-3})$, exist for all $s$.  Moreover, since our assumptions on $\rho_n$ make $G_{n+1}(s,x,y,z)$ $C^2$ in all arguments, the solution, 
{\fontsize{8}{7}\begin{equation}
X_{n+1}(s;s_0,X_0)=(x_{n+1}(s;s_0,X_0),y_{n+1}(s;s_0,X_0),z_{n+1}(s;s_0,X_0)),
\end{equation}}
of $(\ref{xyzchar1})-(\ref{xyzchar3})$, subject to $(\ref{geninitcondchar1-3})$ is, in fact, $C^2$ with respect to all arguments (See theorem 2.10, p.46 of \cite{teschl}).  In particular, $x_{n+1}(s;s_0,X_0),y_{n+1}(s;s_0,X_0)$ and $z_{n+1}(s;s_0,X_0)$ are differentiable with respect to the initial conditions, $s_0$ and $X_0$.    
\end{proof}

Note, that the solution $X_{n+1}(s,X_0)=(x_{n+1}(s,X_0),y_{n+1}(s,X_0),z_{n+1}(s,X_0))$ of $(\ref{xyzchar1})-(\ref{xyzchar3})$, subject to $(\ref{specinitcondchar1-3}),$ satisfies $x(s)=s$, and \begin{eqnarray}
 \label{zsol}z_{n+1}(s;\theta_0)= \rho_{0}(\theta_{0}) e^{k\int_{0}^{s}\int_{0}^{2\pi}\rho_{n}(\tau,\hat{\theta})cos(\hat{\theta}-y_{n+1}(\tau;\theta_0))d\hat{\theta}\mbox{ }d\tau} .\end{eqnarray}  In particular, since $\rho_0$ is nonnegative, so is $z_{n+1}(s;\theta_0)$.


By the existence and uniqueness established in Lemma 1, if we define $\phi_n:\mathbb{R}^2\rightarrow\mathbb{R}^2$ by $\phi_n(s,\theta_0)=(s,y_{n+1}(s;\theta_0))$ then $\phi_n$ is onto and invertible with $\phi_n^{-1}(s,\theta)=(s,y_{n+1}(0;s,\theta)).$ In fact, we have the following corollary: 

\newtheorem*{cor1}{Corollary 1}
\begin{cor1}
	If $\rho_0(\theta)$ and $\rho_n(t,\theta)$ are $C^2$ in all arguments then there exists a $C^2$ solution of (\ref{appxseq}) together with (\ref{IC_z}) on $[0,T]\times\mathbb{R}$.  
\end{cor1}
\begin{proof}
As we have just noted, the mapping $\phi_n:\mathbb{R}^2\rightarrow{\mathbb{R}^2}$ which carries $(s,\theta_0)$ to $(s,y_{n+1}(s;\theta_0))$, where $y_{n+1}(s;\theta_0)$ is the unique solution of (\ref{xyzchar2}) subject to $y_{n+1}(0)=\theta_0$, is both onto $\mathbb{R}^2$ and one-to-one, and hence invertible.  Moreover, by Lemma 1, $\phi_n$ is $C^2$.  In fact, $\frac{\partial\phi_n}{\partial s}(s,\theta_0)=(1,[G_{n+1}]_2(s,y_{n+1}(s;\theta_0)),$ and $\frac{\partial\phi_n}{\partial \theta_0}(s,\theta_0)=(0,\mathrm{exp}\left(\int_0^s{\frac{\partial [G_{n+1}]_2}{\partial y_{n+1}}(\nu,y_{n+1}(\nu;\theta_0))d\nu}\right).$  Hence, the Jacobian of $\phi_n$,
{\fontsize{8}{7}
\begin{eqnarray}
\begin{vmatrix}
                1 & 0\\ [G_{n+1}]_2(s,y_{n+1}(s;\theta_0)) & \mathrm{exp}\left(\int_0^s{\frac{\partial [G_{n+1}]_2}{\partial y_{n+1}}(\nu,y_{n+1}(\nu;\theta_0))d\nu}\right)
\end{vmatrix} \nonumber= \\
\mathrm{exp}\left(\int_0^s{\frac{\partial [G_{n+1}]_2}{\partial y_{n+1}}(\nu,y_{n+1}(\nu;\theta_0))d\nu}\right),
\end{eqnarray}
}
is nonzero (since $\frac{\partial [G_{n+1}]_2}{\partial y_{n+1}}(\nu,y_{n+1}v(\nu;\theta_0))$ is bounded for $\nu\in[0,s]$).  Hence, by the inverse function theorem \cite{evans}, $\phi_n^{-1}$ is $C^2$ at each point in $\mathbb{R}^2$.  Moreover, we can calculate the total derivative of $\phi_n^{-1}$ at $(s,\theta)$ as
\begin{eqnarray}D\phi_n^{-1}=&&\frac{1}{\mathrm{exp}\left(\int_0^s{\frac{\partial [G_{n+1}]_2}{\partial y_{n+1}}(\nu,y_{n+1}(\nu;\theta_0))d\nu}\right)}\nonumber\\
&&\times\begin{bmatrix}\mathrm{exp}\left(\int_0^s{\frac{\partial [G_{n+1}]_2}{\partial y_{n+1}}(\nu,y_{n+1}(\nu;\theta_0))d\nu}\right) & 0\\
				-[G_{n+1}]_2(s,y_{n+1}(s;\theta_0)) & 1\end{bmatrix} \label{jacphi} 
\end{eqnarray}
   where $\theta_0=y_{n+1}(0;s,\theta).$                          
 Using (\ref{jacphi}), we see that 
 \begin{equation}
     \rho_{n+1}(s,\theta):=z_{n+1}(\phi_n^{-1}(s,\theta))\label{rhoSolution}
 \end{equation} is a $C^2$ global solution of (\ref{appxseq}). Moreover, $\rho_{n+1}(0,\theta_0)=z_{n+1}(\phi_n^{-1}(0,\theta_0))=z_{n+1}(0;\theta_0)=\rho_0(\theta_0).$ 
 \end{proof}
 
 Note also, the proof of Corollary 1 shows that $D\phi_n^{-1}$ is uniformly bounded on $[0,T]\times\mathbb{R}$ since since $\rho_n$ is bounded on $[0,2\pi].$ In fact, if $\rho_n$ is normalized on $[0,2\pi],$ 
 \begin{equation}
 \left|\frac{1}{\mathrm{exp}\left(\int_0^s{\frac{\partial [G_{n+1}]_2}{\partial y_{n+1}}(\nu,y_{n+1}(\nu;\theta_0))d\nu}\right)}\right|\leq e^{kT}.
 \end{equation}
 Hence, $\phi_n^{-1}$ is uniformly continuous. 
 
 \newtheorem*{cor2}{Corollary 2}
\begin{cor2}
	If $\rho_n$ is a $C^2$ solution of (\ref{appxseq}) together with (\ref{IC_z}) on $[0,T]\times\mathbb{R}$, then $\frac{\partial \rho_{n+1}}{\partial t}(\phi_n(t,\theta_0))$ and $\frac{\partial \rho_{n+1}}{\partial \theta}(\phi_n(t,\theta_0))$ satisfy the characteristic equations (\ref{p1p2char1}) and (\ref{p1p2char2}) together with the initial data 
{\fontsize{8}{7}
\begin{eqnarray}
\label{ICp1}p_{n+1}(0;\theta_0)=&&-\frac{\partial \rho_0}{\partial \theta}(\theta_0)\left(k \int_0^{2\pi}{\rho_0(\hat{\theta})\sin(\hat{\theta}-\theta_0)\mbox{ }d\hat{\theta}}\right)\nonumber\\
&&+\rho_0(\theta_0)\left(k \int_0^{2\pi}{\rho_0(\hat{\theta})\cos(\hat{\theta}-\theta_0)\mbox{ }d\hat{\theta}}\right),\\
\label{ICp2}q_{n+1}(0;\theta_0)=&&\rho_0'(\theta_0),
\end{eqnarray}
}
respectively. That is $\frac{\partial \rho_{n+1}}{\partial t}(t,\theta)=p_{n+1}(\phi_n^{-1}(t,\theta))$ and $\frac{\partial \rho_{n+1}}{\partial \theta}(t,\theta)=q_{n+1}(\phi_n^{-1}(t,\theta))$
\label{cor2}\end{cor2}
\begin{proof}See Theorem 1, page 99 of \cite{evans}.\end{proof}
				


				
\newtheorem*{lem3}{Lemma 3}
\begin{lem3}\label{Lemma3}
If $\rho_n$ and $\rho_0$ are nonnegative and $2\pi$-periodic in $\theta$, and if $  \int_{0}^{2\pi}\rho_{n}(t,\theta) d\theta=1,$ then 
\begin{itemize}
\item[i]  $ \rho_{n+1}(t,\theta)\geq{0}$.\\
\item[ii] $ \rho_{n+1}(t,\theta) $ is  $ 2\pi $-periodic in $ \theta $.\\
\item[iii] $ \int_{0}^{2\pi}\rho_{n+1}(t,\theta) d\theta=1.$
\end{itemize}
\end{lem3}
\begin{proof}
\begin{itemize}\begin{flushleft}
\item[i.] Since $ \rho_{n+1}(t,\theta) = z_{n+1}(\phi_n^{-1}(t,\theta)), $ where $z_{n+1}$ denotes the third component of the solution of $(\ref{xyzchar1})-(\ref{xyzchar3})$ under the initial data  $ (0, y_{n+1}(0), z_{n+1}(0)) = (0, \theta_{0}, \rho_0(\theta_0)),$ by (\ref{zsol}), $\rho_{n+1}(t,\theta)\geq{0}$, provided $\rho_0(\theta)\geq{0}.$
\item[ii.] Given any $t$ and $\theta$, we evaluate $ \rho_{n+1}(t,\theta) $ and $ \rho_{n+1}(t,\theta+2\pi) $ using characteristics. First note that if $(s, y_{n+1}(s;\theta_0))$ is the unique characteristic curve through $(t,\theta)$, i.e. if $y_{n+1}(t;\theta_0)=\theta$, then $(s,y_{n+1}(s;\theta_0)+2\pi)$ is the unique characteristic curve through $(t,\theta+2\pi)$. That is, $(s,y_{n+1}(s;\theta_0+2\pi))=(s,y_{n+1}(s;\theta_0)+2\pi).$ Indeed,
{\fontsize{8}{7}
\begin{align}
    \frac{d}{ds}(y_{n+1}(s;\theta_0+2\pi)) &= k\int_{0}^{2\pi}\rho_n(s,\hat{\theta})\sin(\hat{\theta}-y_{n+1}(s;\theta_0+2\pi))d\hat{\theta}\\
    \frac{d}{ds}(y_{n+1}(s;\theta_0)+2\pi) &= k\int_{0}^{2\pi}\rho_n(s,\hat{\theta})\sin(\hat{\theta}-y_{n+1}(s;\theta_0))d\hat{\theta}\\
    &=k\int_{0}^{2\pi}\rho_n(s,\hat{\theta})\sin(\hat{\theta}-y_{n+1}(s;\theta_0)+2\pi)d\hat{\theta}
\end{align}}
That is $y_{n+1}(s;\theta_0+2\pi)$ and $y_{n+1}(s;\theta_0)+2\pi$ solve the same initial value problem. By uniqueness, \begin{equation}
y_{n+1}(s;\theta_0+2\pi)=y_{n+1}(s;\theta_0)+2\pi.\end{equation} 
In other words, $y_{n+1}(s;\theta_0)  $ $mod$ $ 2\pi$  is $2\pi$-periodic (in $\theta$).
Since cosine and $\rho_0$ are $2\pi$-periodic in $\theta$. From this observation and (\ref{zsol}), we see that $z_{n+1}(s;\theta_0)=z_{n+1}(s;\theta_0+2\pi)$ for all $s$.  Finally, given $(t,\theta)\in[0,T]\times\mathbb{R}$ there exists $\theta_0$ so that $y(t;\theta_0)=\theta$. Hence,
\begin{align}\rho_{n+1}(t,\theta)&=z_{n+1}(t;\theta_0)\\
&=z_{n+1}(t;\theta_0+2\pi)\\
&=\rho_{n+1}(t,y_{n+1}(t;\theta_0+2\pi))\\
&=\rho_{n+1}(t,y_{n+1}(t;\theta_0)+2\pi)\\
&=\rho_{n+1}(t,\theta+2\pi)
\end{align}
\item[iii.] Recall (\ref{appxseq}) can be expressed as 
{\fontsize{8}{7}\begin{equation}
\frac{\partial \rho_{n+1}}{\partial t} (t,\theta)= -\frac{\partial}{\partial \theta}\left(k \rho_{n+1}(t,
\theta)\int_0^{2\pi}{\rho_n(t,\hat{\theta})\sin(\hat{\theta}-\theta)d\hat{\theta}}\right).
\end{equation}}
Integrating both sides with respect to $ \theta $ from 0 to $ 2\pi $, we get 
{\fontsize{7}{6}\begin{equation}
\frac{\partial }{\partial t} \int_{0}^{2\pi}\rho_{n+1}(t,\theta)d\theta= -\left(k \rho_{n+1}(t,
\theta)\int_0^{2\pi}{\rho_n(t,\hat{\theta})\sin(\hat{\theta}-\theta)d\hat{\theta}}\right)_{\theta=0}^{\theta=2\pi}=0.
\end{equation}}
Thus $ \int_{0}^{2\pi}\rho_{n+1}(t,\theta)d\theta $ is constant with respect to $ t $, and so by our assumptions on $\rho_n$, (iii) holds. 
	\end{flushleft}
\end{itemize}  

\end{proof}

\newtheorem*{cor3}{Corollary 3}
\begin{cor3}
	If $\rho_n(t,\theta)$ is a continuous then a $C^1$ solution of (\ref{appxseq}) together with (\ref{IC_z}) on $[0,T]\times\mathbb{R}$ is unique.  
\end{cor3}

\begin{proof}
 Suppose $u_1(t,\theta)$ and $u_2(t,\theta)$ are $C^2$, solutions of (\ref{appxseq}) and (\ref{IC_z}) on $[0,T]\times\mathbb{R}$, where $\rho_n(t,\theta)$ is $2\pi$-periodic in $\theta$. Then, by Lemma 3, $u_1(t,\theta)$ and $u_2(t,\theta)$ are $2\pi$-periodic in $\theta$. Consider the integral of the difference of these solutions squared with respect to $\theta$. Differentiating with respect to $t$ we find
{\fontsize{7}{6}
\begin{align}
    \frac{d}{dt}\int_{0}^{2\pi}&(u_1-u_2)^2(t,\theta)\ d\theta=\int_{0}^{2\pi}2(u_1-u_2)(t,\theta)\frac{\partial(u_1-u_2)}{\partial t}(t,\theta)\ d\theta\nonumber\\
    =&
    -2k\int_{0}^{2\pi}(u_1-u_2)\frac{\partial (u_1-u_2)}{\partial\theta}(t,\theta)\int_{0}^{2\pi}\rho_n(t,\hat{\theta})\sin(\hat{\theta}-\theta)d\hat{\theta}\ d\theta
    \nonumber\\&+2k\int_{0}^{2\pi}(u_1-u_2)^2(t,\theta)\int_{0}^{2\pi}\rho_n(t,\hat{\theta})\cos(\hat{\theta}-\theta)\ d\hat{\theta}\ d\theta\nonumber\\
    =\label{integral_by_parts}&-k\int_{0}^{2\pi}\frac{\partial(u_1-u_2)^2}{\partial\theta}(t,\theta)\int_{0}^{2\pi}\rho_n(t,\hat{\theta})\sin(\hat{\theta}-\theta)d\hat{\theta}d\theta
    \\&+2k\int_{0}^{2\pi}(u_1-u_2)^2(t,\theta)\int_{0}^{2\pi}\rho_n(t,\hat{\theta})\cos(\hat{\theta}-\theta)d\hat{\theta}d\theta.
\end{align}
}
Then, through use of integration by parts on (\ref{integral_by_parts})
{
\begin{align}
\mu&=k\int_{0}^{2\pi}\rho_n(t,\hat{\theta})\sin(\hat{\theta}-\theta)d\hat{\theta}\\
d\mu &=-k\int_{0}^{2\pi}\rho_n(t,\hat{\theta})\cos(\hat{\theta}-\theta)d\hat{\theta}d\theta\\
d\lambda &= -\frac{\partial(u_1-u_2)^2(t,\theta)}{\partial\theta}d\theta \\
\lambda&=-(u_1-u_2)^2(t,\theta)
\end{align}
}
Since sine, $u_1$, and $u_2$ are $2\pi$-periodic in $\theta$, $\mu\lambda\big|_0^{2\pi} =0$, and (\ref{integral_by_parts}) reduces to 
\begin{equation}
-k\int_{0}^{2\pi}(u_1-u_2)^2(t,\theta)\int_{0}^{2\pi}\rho_n(t,\hat{\theta})\cos(\hat{\theta}-\theta)\ d\hat{\theta}\ d\theta.
\end{equation}
Thus,
{\fontsize{7}{6}
\begin{align}
    \frac{d}{dt}\int_{0}^{2\pi}(u_1-u_2)^2(t,\theta)\ d\theta =&k\int_{0}^{2\pi}(u_1-u_2)^2(t,\theta)\int_{0}^{2\pi}\rho_n(t,\hat{\theta})\cos(\hat{\theta}-\theta)\ d\hat{\theta}\ d\theta\nonumber\\
    \leq& 2k\pi\|\rho_n\|_{\infty} \int_{0}^{2\pi}(u_1-u_2)^2(t,\theta)\ d\theta,
\end{align}}
where $\left\|\cdot\right\|_{\infty}$ denotes the supremum on $[0,T]\times[0,2\pi],$ and therefore
\begin{equation}
\int_{0}^{2\pi}(u_1-u_2)^2(t,\theta)\ d\theta \equiv \int_{0}^{2\pi}(u_1-u_2)^2(0,\theta)\ d\theta =0.
\end{equation}
It follows that $C^1$ solutions to (\ref{appxseq}) and (\ref{IC_z}) on $[0,T]\times\mathbb{R}$ are unique.
\end{proof}

Now we consider the existence and uniqueness of solutions $ p_{n+1} $ and $ q_{n+1} $ to the  characteristic equations and initial conditions associated with the first derivatives of $\rho_{n+1}.$ Since we have seen that $x(s;\theta_0)=s$ we replace (\ref{p1p2char1}) and (\ref{p1p2char2}) with
{\fontsize{8}{7}
\begin{eqnarray}
\label{p1p2char1_fin}\dot{p}_{n+1}(s;\theta_0)=&& -q_{n+1}(s;\theta_0)k\int_0^{2\pi}{\frac{\partial}{\partial s}\rho_{n}(s,\hat{\theta})\sin(\hat{\theta}-y_{n+1}(s;\theta_0))\ d\hat{\theta}}\nonumber,\\
&&+ z_{n+1}(s;\theta_0)k\int_0^{2\pi}{\frac{\partial}{\partial s}\rho_{n}(s,\hat{\theta})\cos(\hat{\theta}-y_{n+1}(s;\theta_0))\ d\hat{\theta}}\\
&&+ p_{n+1}(s;\theta_0)k\int_0^{2\pi}{\rho_{n}(s,\hat{\theta})\cos(\hat{\theta}-y_{n+1}(s;\theta_0))\ d\hat{\theta}}\nonumber,\\
\label{p1p2char2_fin}\dot{q}_{n+1}(s;\theta_0)=&& 2q_{n+1}(s;\theta_0)k\int_0^{2\pi}{\rho_{n}(s,\hat{\theta})\cos(\hat{\theta}-y_{n+1}(s;\theta_0))\ d\hat{\theta}}\nonumber\\
&&+z_{n+1}(s)k\int_0^{2\pi}{\rho_{n}(s,\hat{\theta})\sin(\hat{\theta}-y_{n+1}(s;\theta_0))\ d\hat{\theta}},
\end{eqnarray}
}
which will be solved subject to (\ref{ICp1}) and (\ref{ICp2}). 
\newtheorem*{lem4}{Lemma 4}
\begin{lem4}
If $\rho_0(\theta)$ and $\rho_n(t,\theta)$ are $C^2$ in all arguments, there exists a unique, continuously differentiable solution to (\ref{p1p2char1_fin})-(\ref{p1p2char2_fin}) subject to initial conditions (\ref{ICp1})-(\ref{ICp2}) on $[0,T]\times\mathbb{R}$. 
\end{lem4}
\begin{proof}
 Let $[G_{n+1}]_4(s,p,q)$ and $[G_{n+1}]_5(s,p,q)$ be defined by the right-hand side of (\ref{p1p2char1_fin}) and (\ref{p1p2char2_fin}), respectively.  Since $x_{n+1}(s;s_0,X_0)$, $y_{n+1}(s,s_0,X_0)$, and $z_{n+1}(s;s_0,X_0)$ are $C^2$ by Lemma 1, and $\rho_{n}$ is $C^2$ in each of its arguments by hypothesis, we see that $[G_{n+1}]_5$ is $C^2$ in all variables by the Leibniz Rule.  To see that $[G_{n+1}]_4$ is $C^2$, despite its dependence on $\frac{\partial\rho_{n}}{\partial s}$, note that in case $n=0$, $\rho_{n}(x,\hat{\theta})=\rho_0(\hat{\theta)}$, so $\frac{\partial \rho_{n}}{\partial s}\equiv0$, and $[G_{n+1}]_4(s,p_{n+1},q_{n+1})$ reduces to an expression that is clearly $C^2$ is all variables. In case $n\geq1$, 
{\fontsize{7}{6}
\begin{eqnarray}\label{byparts}
\frac{\partial}{\partial s}\rho_{n}(s,\hat{\theta})&=&-\frac{\partial}{\partial \hat{\theta}}\left(\rho_{n}(s,\hat{\theta})k\int_0^{2\pi}{\rho_{n-1}(s,\hat{\hat{\theta}})\sin(\hat{\hat{\theta}}-\hat{\theta})d\hat{\hat{\theta}}}\right).
\end{eqnarray}}
Substituting the right hand side of (\ref{byparts}) into (\ref{p1p2char1_fin}) and integrating by parts we have,
{\fontsize{6}{5}
\begin{eqnarray}\label{G4}
\dot{p}_{n+1}(s;\theta_0)=&&-q_{n+1}(s;\theta_0)k\int_0^{2\pi}{\rho_{n}(s,\hat{\theta})\cos(\hat{\theta}-y_{n+1}(s;\theta_0))\int_0^{2\pi}{\rho_{n-1}(s,\hat{\hat{\theta}})\sin(\hat{\hat{\theta}}-\hat{\theta})\ d\hat{\hat{\theta}}}\ d\hat{\theta}}\nonumber\\
&&-z_{n+1}(s;\theta_0)k\int_0^{2\pi}{\rho_{n}(s,\hat{\theta})\sin(\hat{\theta}-y_{n+1}(s;\theta_0))\int_0^{2\pi}{\rho_{n-1}(s,\hat{\hat{\theta}})\sin(\hat{\hat{\theta}}-\hat{\theta})d\hat{\hat{\theta}}}\ d\hat{\theta}}\nonumber\\
&&+ p_{n+1}(s;\theta_0)k\int_0^{2\pi}{\rho_{n}(s,\hat{\theta})\cos(\hat{\theta}-y_{n+1}(s;\theta_0))\ d\hat{\theta}}.
\end{eqnarray}
}
So that $[G_{n+1}]_4$ is also $C^2$ in $p, q$ and $s$.  In particular, $ [G_{n+1}]_4 $ and $ [G_{n+1}]_5 $ are linear in $p$ and $q$ and hence uniformly Lipschitz continuous on $[0,T]\times\mathbb{R}^2$ with respect to $p$ and $q$. It follows that there exists a unique global solution of (\ref{p1p2char1_fin})-(\ref{p1p2char2_fin}) subject to (\ref{ICp1})-(\ref{ICp2}) (See corollary 2.6, p.41 of \cite{teschl}).  Moreover a solution of (\ref{p1p2char1_fin})-(\ref{p1p2char2_fin}) is $C^2$ with respect to both $s$ and its initial data, $p_{n+1}(0), q_{n+1}(0)$.  Since this data is a continuously differentiable function of $\rho_0(\theta_0)$ and $\frac{\partial\rho_0}{\partial \theta}(\theta_0),$ which are continuously differentiable with respect to $\theta_0,$ we see that the solution of (\ref{p1p2char1_fin})-(\ref{p1p2char2_fin}) subject to (\ref{ICp1})-(\ref{ICp2}) is continuously differentiable with respect to $s$ and $\theta_0$.  
\end{proof}

Alternatively, we can give an explicit formula for the solution of (\ref{p1p2char1_fin})-(\ref{p1p2char2_fin}) subject to (\ref{ICp1})-(\ref{ICp2}).  In particular, the solution of (\ref{p1p2char2_fin}) subject to (\ref{ICp2}) is
{\fontsize{7}{6}
\begin{eqnarray}
q_{n+1}(t;\theta_0)=&&\mathrm{e}^{\int_0^t{2f_{n+1}(s;\theta_0)ds}}\nonumber\\
&&\times\left[\int_0^t{g_{n+1}(s;\theta_0)z_{n+1}(s;\theta_0)\mathrm{e}^{-\int_0^s{2f_{n+1}(\nu;\theta_0)d\nu}}ds}+\rho_0'(\theta_0)\right]\\
=&&z_{n+1}(t;\theta_0)^{2}\left[\int_0^t{\frac{g_{n+1}(s;\theta_0)}{z_{n+1}(s;\theta_0)}ds}+\frac{\rho_0'(\theta_0)}{\rho_0(\theta_0)^{2}}\right],
\label{p1p2char2sol}
\end{eqnarray}}
provided $z_{n+1}(0;\theta_0)=\rho_0(\theta_0)\neq{0}$, where \begin{eqnarray}
f_{n+1}(s;\theta_0):&=&k\int_0^{2\pi}{\rho_{n}(s,\hat{\theta})\cos(\hat{\theta}-y_{n+1}(s;\theta_0))d\hat{\theta}}\nonumber\\
&=&\frac{\partial [G_{n+1}]_2}{\partial y}(s,y_{n+1}(s,\theta))
\end{eqnarray}and 
\begin{eqnarray}
g_{n+1}(s;\theta_0):&=&k\int_0^{2\pi}{\rho_{n}(s,\hat{\theta})\sin(\hat{\theta}-y_{n+1}(s;\theta_0))d\hat{\theta}}\nonumber
\\&=&[G_{n+1}]_2(s,y_{n+1}(s,\theta)),
\end{eqnarray}
and, in case $z_{n+1}(0;\theta_0)=\rho_0(\theta_0)=0$,
\begin{eqnarray}
q_{n+1}(t;\theta_0)&=&\rho_0'(\theta_0)\mathrm{e}^{\int_0^t{2f_{n+1}(s;\theta_0)ds}}.
\label{p2sol2}
\end{eqnarray}
Hence,
{\fontsize{8}{7}
\begin{eqnarray}
\label{p1b1} p_{n+1}(t;\theta_0)=&&\mathrm{e}^{\int_0^t{f_{n+1}(s;\theta_0)ds}}\nonumber\\
&&\times\left[\int_0^t{c_{n+1}(s;\theta_0)\mathrm{e}^{-\int_0^s{f_{n+1}(\nu;\theta_0)d\nu}}ds}+p_{n+1}(0;\theta_0)\right]\label{p1p2char1solgen}\\
\label{p1b2}=&&z_{n+1}(t;\theta_0)\left[\int_0^t{\frac{c_{n+1}(s;\theta_0)}{z_{n+1}(s;\theta_0)}ds}+\frac{p_{n+1}(0;\theta_0)}{\rho_0(\theta_0)}\right].
\label{p1p2char1sol}
\end{eqnarray}}
where $(\ref{p1p2char1sol})$ holds provided $z_{n+1}(0;\theta_0)=\rho_0(\theta_0)\neq{0},$ and 
{\fontsize{6}{5}
\begin{align}
    c_{n+1}(s;\theta_0)=&-q_{n+1}(s;\theta_0)k\int_0^{2\pi}{\rho_{n}(s;\hat{\theta})\cos(\hat{\theta}-y_{n+1}(s;\theta_0))\int_0^{2\pi}{\rho_{n-1}(s,\hat{\hat{\theta}})\sin(\hat{\hat{\theta}}-\hat{\theta})d\hat{\hat{\theta}}}d\hat{\theta}}\nonumber\\
&-z_{n+1}(s;\theta_0)k\int_0^{2\pi}{\rho_{n}(s,\hat{\theta})\sin(\hat{\theta}-y_{n+1}(s;\theta_0))\int_0^{2\pi}{\rho_{n-1}(s,\hat{\hat{\theta}})\sin(\hat{\hat{\theta}}-\hat{\theta})d\hat{\hat{\theta}}}d\hat{\theta}}
\end{align}
}

Now, from (\ref{zsol}) and Lemma 2 parts (i) and (iii) we see that for all $\theta_0,$ 
{\fontsize{8}{7}
\begin{eqnarray}
\label{b1}|z_{n+1}(t;\theta_0)|&\leq& \rho_0(\theta_0)\mathrm{e}^{kt}\\
\label{lb1}|z_{n+1}(t;\theta_0)|&\geq& \rho_0(\theta_0)\mathrm{e}^{-kt}\\
\label{b2}|f_{n+1}(t;\theta_0)|&\leq& k\\
\label{b3}|g_{n+1}(t;\theta_0)|&\leq& k\\
|q_{n+1}(t;\theta_0)|&\leq&|z_{n+1}(t;\theta_0)^{2}|\left(\left|\int_0^t{\frac{g_{n+1}(s;\theta_0)}{z_{n+1}(s;\theta_0)}ds}\right|+\left|\frac{\rho_0'(\theta_0)}{\rho_0(\theta_0)^{2}}\right|\right)\nonumber\\
&=&e^{\int_{0}^{t}2f_{n+1}(s; \theta_0)\mbox{ }ds}\left(\rho_0(\theta_0)\int_0^t{\left|\frac{g_{n+1}(s;\theta_0)}{e^{\int_{0}^{s}f_{n+1}(\tau; \theta_0)\mbox{ }d\tau}}\right|ds}+\left|\rho_0'(\theta_0)\right|\right)\nonumber\\
&\leq&e^{2kt}\left(\rho_0(\theta_0)e^{kt}+\left|\rho_0'(\theta_0)\right|\right)\quad z_{n+1}(0)\neq(0)\label{bq1}\\
\label{bq2}|q_{n+1}(t;\theta_0)|&\leq&{\rho_0'(\theta_0)\mathrm{e}^{2kt}},\quad z_{n+1}(0)=(0).
\end{eqnarray}}
 Note that since $\rho_0$ is $C^2$ and $2\pi$-periodic, (\ref{b1})-(\ref{b3}) provide a uniform bound, $B_z(T)$ for $\rho_{n+1}(t,\theta)=z_{n+1}(\phi_n^{-1}(t,\theta))$ on $[0,T]\times \mathbb{R}$.  Also, by Corollary 2, (\ref{bq1})-(\ref{bq2}) provide a uniform bound $B_q(T)$ for $\frac{\partial \rho_{n+1}}{\partial \theta}(t,\theta)=q_{n+1}(\phi_n^{-1}(t,\theta))$ on $[0,T]\times \mathbb{R}$. That is, $\left\{\rho_n\right\}_{n=1}^{\infty}$ and $\left\{\frac{\partial \rho_n}{\partial \theta}\right\}_{n=1}^{\infty}$ are uniformly bounded. It follows from (\ref{appxseq}) that 
\begin{eqnarray}\label{b5}
\left|\frac{\partial \rho_n}{\partial t}(t,\theta)\right|&\leq&k(B_{q}(T)+B_z(T)),
\end{eqnarray}
for all $n\in \mathbb{N}$ and $(t,\theta)\in[0,T]\times \mathbb{R}.$ Thus we have the following lemma:

\newtheorem*{lem5}{Lemma 5}
\begin{lem5}
The sequence $\left\{\rho_n\right\}_{n=0}^{\infty}$ is uniformly bounded and equicontinuous on $[0,T]\times \mathbb{R}$.  
\end{lem5}

In establishing a global solution of (\ref{PDE}), we would like to show also that $\frac{\partial \rho_n}{\partial t }$ and $\frac{\partial \rho_n}{\partial \theta}$ are uniformly bounded and equicontinuous on ${[0,T]}\times\mathbb{R}$. Toward this goal, we have already established uniform bounds. For equicontinuity, it suffices to show $D\frac{\partial \rho_n}{\partial t }=\left[\frac{\partial^2 \rho_n}{\partial t^2 }, \frac{\partial^2 \rho_n}{\partial t \partial \theta }\right]$ and $D\frac{\partial \rho_n}{\partial \theta }=\left[\frac{\partial^2 \rho_n}{\partial t \partial \theta }, \frac{\partial^2 \rho_n}{\partial \theta^2 }\right]$ are uniformly bounded on $[0,T]\times \mathbb{R}$ for all $T>0$. We have established the existence of a $C^2$ solution $\rho_n(t,\theta)$ of (\ref{appxseq}) on $[0,T]\times\mathbb{R}$. Moreover, along the projected characteristics (\ref{xyzchar1})-(\ref{xyzchar2}), $\rho_n$ and its partial derivatives $\frac{\partial \rho_n}{\partial t }$ and $\frac{\partial \rho_n}{\partial \theta}$, satisfy (\ref{xyzchar3}), (\ref{p1p2char1_fin}), and (\ref{p1p2char2_fin}), together with (\ref{IC_z}), (\ref{ICp1}), and (\ref{ICp2}), respectively. That is, 
\begin{align}
    \label{bounded1}\frac{\partial \rho_{n+1}}{\partial t}(t,y_{n+1}(t;\theta_0)) &= p_{n+1}(t;\theta_0)\\
    \label{bounded2}\frac{\partial \rho_{n+1}}{\partial \theta}(t,y_{n+1}(t;\theta_0)) &= q_{n+1}(t;\theta_0)
\end{align}
Thus,
\begin{widetext}
{\begin{eqnarray}\label{totalDerivativep1}
    D\left[\frac{\partial \rho_n}{\partial t }(t,y_{n+1}(t;\theta_0))\right]
    &=\left[\frac{\partial^2 \rho_n}{\partial t^2 }(t,y_{n+1}(t;\theta_0)), \frac{\partial^2 \rho_n}{\partial t \partial \theta }(t,y_{n+1}(t;\theta_0))\right]\begin{bmatrix}
                1 & 0\\ G_2(s,y_{n+1}(s;\theta_0) & \mathrm{exp}\left(\int_0^s{\frac{\partial G_2}{\partial y_{n+1}}(\nu,y_{n+1}(\nu;\theta_0))d\nu}\right)
    \end{bmatrix}\nonumber\\
    &=\left[\dot{p}_{n+1}(t;\theta_0),\frac{\partial p_{n+1}}{\partial \theta_0}(t;\theta_0)\right],
\end{eqnarray}}
which yields,\begin{eqnarray}
    \left[\frac{\partial^2 \rho_n}{\partial t^2 }(t,y_{n+1}(t;\theta_0)), \frac{\partial^2 \rho_n}{\partial t \partial \theta }(t,y_{n+1}(t;\theta_0))\right] =&&
   \frac{1}{\mathrm{exp}\left(\int_0^s{\frac{\partial G_2}{\partial y_{n+1}}(\nu,y_{n+1}(\nu;\theta_0))d\nu}\right)}\nonumber\\
   &&\times\left[\dot{p}_{n+1}(t;\theta_0),\frac{\partial p_{n+1}}{\partial\theta_0}(t;\theta_0)\right]\begin{bmatrix}\mathrm{exp}\left(\int_0^s{\frac{\partial G_2}{\partial y_{n+1}}(\nu,y_{n+1}(\nu;\theta_0))d\nu}\right) & 0\\
				-G_2(s,y_{n+1}(s;\theta_0)) & 1\end{bmatrix}
\end{eqnarray}
Also,
\begin{eqnarray}
    D\left[\frac{\partial \rho_n}{\partial \theta }(t,y_{n+1}(t;\theta_0))\right]
    =&&\left[\frac{\partial \rho_n}{\partial t\partial\theta }(t,y_{n+1}(t;\theta_0)), \frac{\partial^2 \rho_n}{\partial \theta^2 }(t,y_{n+1}(t;\theta_0))\right]\begin{bmatrix}
                1 & 0\\ G_2(s,y_{n+1}(s;\theta_0) & \mathrm{exp}\left(\int_0^s{\frac{\partial G_2}{\partial y_{n+1}}(\nu,y_{n+1}(\nu;\theta_0))d\nu}\right)
    \end{bmatrix}\nonumber\\
    =&&\left[\dot{q}_{n+1}(t;\theta_0), \frac{\partial q_{n+1}}{\partial \theta_0}(t;\theta_0)\right],
\end{eqnarray}
which yields,
\begin{eqnarray}
  \left[\frac{\partial \rho_n}{\partial t\partial\theta }(t,y_{n+1}(t;\theta_0)), \frac{\partial^2 \rho_n}{\partial \theta^2 }(t,y_{n+1}(t;\theta_0))\right]=&&   \frac{1}{\mathrm{exp}\left(\int_0^s{\frac{\partial G_2}{\partial y_{n+1}}(\nu,y_{n+1}(\nu;\theta_0))d\nu}\right)}\nonumber\\
  &&\times\left[\dot{q}_{n+1}(t;\theta_0), \frac{\partial q_{n+1}}{\partial \theta_0}(t;\theta_0)\right]\begin{bmatrix}\mathrm{exp}\left(\int_0^s{\frac{\partial G_2}{\partial y_{n+1}}(\nu,y_{n+1}(\nu;\theta_0))d\nu}\right) & 0\\
				-G_2(s,y_{n+1}(s;\theta_0)) & 1\end{bmatrix}\label{totalDerivativep2}
\end{eqnarray}
\end{widetext}
We will proceed by showing the the total derivative of $\frac{\partial \rho_n}{\partial \theta}$, i.e. (\ref{totalDerivativep2})), is uniformly bounded, independent of $n$, on $[0,T]\times \mathbb{R}$. It will then follow that the sequence $\left\{\frac{\partial \rho_n}{\partial \theta}\right\}_{n=0}^{\infty},$ is equicontinuous. We have already seen that $D\phi_n^{-1}$ is uniformly bounded, so it remains to show that the partial derivatives $\frac{\partial q_{n+1}}{\partial \theta_0}(t;\theta_0)$ and $\dot{q}_{n+1}(t;\theta_0)$ are uniformly bounded for $(t,\theta_0)\in [0,T]\times \mathbb{R}$. 

First we derive some bounds on the derivatives of $y_{n+1}$. If we differentiate (\ref{xyzchar2}) with respect to the initial condition $\theta_0$, we get
{\fontsize{8}{7}\begin{align}
    \frac{\partial \dot{y}_{n+1}}{\partial \theta_0}(s;\theta_0) &= \frac{\partial}{\partial \theta_0}\left(k\int_{0}^{2\pi}\rho_n(s,\hat{\theta})\sin(\hat{\theta}-y_{n+1}(s;\theta_0))d\hat{\theta}\right)\nonumber\\
     &= -\frac{\partial y_{n+1}}{\partial \theta_0}(s;\theta_0)k\int_{0}^{2\pi}\rho_n(s,\hat{\theta})\cos(\hat{\theta}-y_{n+1}(s;\theta_0))d\hat{\theta},
\end{align}}
where we have used $\frac{\partial }{\partial \theta_0}x_{n+1}(t;\theta_0)\equiv0.$ This yields 
{\fontsize{8}{7}\begin{align}\frac{\partial y_n}{\partial\theta_0}(t,\theta_0)&=\mathrm{exp}\left(-\int_0^t f_{n}(s;\theta_0)\ ds\right)\nonumber\\
&=\mathrm{exp}\left(-k\int_0^t{ \int_0^{2\pi}{\rho_{n-1}(s,\hat{\theta})\cos(\hat{\theta}-y_{n}(s;\theta_0))d\hat{\theta}}\ ds}\right),
\end{align}} so
\begin{equation}\label{ParInitYBound}
\left|\frac{\partial y_{n+1}}{\partial \theta_0}(s;\theta_0)\right|\leq e^{kT}.\end{equation}
Recall $f_{n+1}(s;\theta_0)=k\int_0^{2\pi}\rho_n(s,\hat{\theta})\cos(\hat{\theta}-y_{n+1}(s;\theta_0))\ d\hat{\theta},$ and $g_{n+1}(s;\theta_0) = k\int_0^{2\pi}\rho_n(s,\hat{\theta})\sin(\hat{\theta}-y_{n+1}(s;\theta_0))d\hat{\theta}.$ 
Thus, for $s\in[0,T]$ 
$\frac{\partial f}{\partial \theta_0}$ and $\frac{\partial g}{\partial \theta_0}$ are bound in magnitude by $ke^{kT}$ on $[0,T]\times\mathbb{R}$.


Now we look to uniformly bound $\frac{\partial q_{n+1}}{\partial \theta_0}(t; \theta_0)$.  From (\ref{p1p2char2sol}), if $z_{n+1}(0;\theta_0) = \rho_0(\theta_0)  \neq 0$,
{\fontsize{8}{7}\begin{align}
q_{n+1}(t;\theta_0)&= z_{n+1}(t;\theta_0)^{2}\left(\int_0^t{\frac{g_{n+1}(s;\theta_0)}{z_{n+1}(s;\theta_0)}ds}+\frac{\rho_0'(\theta_0)}{\rho_0(\theta_0)^{2}}\right)\\
&=e^{\int_{0}^{t}2f_{n+1}(s; \theta_0)\mbox{ }ds}\left[\rho_0(\theta_0)\int_0^t{\frac{g_{n+1}(s;\theta_0)}{e^{\int_{0}^{s}f_{n+1}(\tau; \theta_0)\mbox{ }d\tau}}}+\rho_0'(\theta_0)\right]
\end{align}}
Since the derivatives of $f_{n+1}$ and $g_{n+1}$ with respect to $\theta_0$ are continuous and bounded in magnitude, since 
\begin{equation}
e^{\int_{0}^{s}f_{n+1}(\tau; \theta_0)\mbox{ }d\tau}\geq e^{-kT},
\end{equation}
and since $\rho_0(\theta_0)$ is $C^2$ and $2\pi$-periodic, we can see from the Leibniz and quotient rules that $\frac{\partial q_{n+1}}{\partial \theta_0} (s;\theta_0)$ is bounded on $[0,T]\times \mathbb{R}$, independent of $n$.
Or, if $z_{n+1}(0;\theta_0) = 0$, then 
\begin{equation}
q_{n+1}(t;\theta_0)=\rho_0'(\theta_0)\mathrm{e}^{\int_0^t{2f_{n+1}(s;\theta_0)ds}},
\end{equation}
so again we see that $\frac{\partial q_{n+1}}{\partial \theta_0} (s;\theta_0)$ is bounded on $[0,T]\times \mathbb{R}$, independent of $n$. 

Also notice from (\ref{p1p2char2}) and the bounds (\ref{b1})-(\ref{bq2}) that $\dot{q}_n(t)$ is uniformly bounded with respect to $t$ on $[0,T]\times \mathbb{R}.$ Hence the sequence $\left\{\frac{\partial \rho_n}{\partial \theta}\right\}_{n=0}^{\infty}$ is equicontinuous on $[0,T]\times \mathbb{R}.$ From (\ref{appxseq}), we see that $\left\{\frac{\partial \rho_n}{\partial t}\right\}_{n=0}^{\infty}$ is also equicontinuous on $[0,T]\times \mathbb{R}.$ Indeed, if $\left\{u_n\right\}_{n=1}^{\infty}$ and $\left\{v_n\right\}_{n=1}^{\infty}$ are uniformly bounded and equicontinuous, then so is $\left\{u_nv_n\right\}_{n=1}^{\infty}$. 

We have established the following lemma.
\newtheorem*{lem6}{Lemma 6}
\begin{lem6}
The sequences $\left\{\frac{\partial \rho_n}{\partial t}\right\}_{n=0}^{\infty}$ and $\left\{\frac{\partial \rho_n}{\partial \theta}\right\}_{n=0}^{\infty}$ are uniformly bounded and equicontinuous on $[0,T]\times \mathbb{R} $.  
\end{lem6}

From Lemmas 5 and 6 and by the Arzela-Ascoli Theorem, we can choose a subsequence $\{n_k\}$ so that $\left\{\rho_{n_k}\right\}_{k=0}^{\infty},$ $\left\{\frac{\partial \rho_{n_k}}{\partial \theta}\right\}_{k=0}^{\infty},$ and $\left\{\frac{\partial \rho_{n_k}}{\partial t}\right\}_{k=0}^{\infty}$ converge uniformly on $[0,T]\times[0,2\pi]$. Then, since $\rho_n(t,\theta)$ is $2\pi$-periodic in $\theta$, convergence is also uniform on $[0,T]\times\mathbb{R}.$ This yields,
\begin{eqnarray}
\lim\limits_{n\to\infty} \rho_{n_k}(t,\theta) &=& \rho(t,\theta)\\
    \lim\limits_{n\to\infty} \frac{\partial}{\partial t} \rho_{n_k}(t,\theta) &=& \frac{\partial}{\partial t}\lim\limits_{n\to\infty} \rho_{n_k}(t,\theta) = \frac{\partial \rho}{\partial t} (t,\theta)\\
    \lim\limits_{n\to\infty} \frac{\partial}{\partial \theta} \rho_{n_k}(t,\theta)  &=& \frac{\partial}{\partial \theta}\lim\limits_{n\to\infty} \rho_{n_k}(t,\theta) = \frac{\partial \rho}{\partial \theta} (t,\theta)
\end{eqnarray}
for $(t,\theta)\in[0,T]\times\mathbb{R}.$
 Furthermore, uniform convergence enables us to pull limits through the integral in (\ref{appxseq}), so the subsequence converges uniformly to a $C^1$ solution of (\ref{PDE}) on $[0,T]\times\mathbb{R}$. Since $T>0$ was arbitrary, we have established the following theorem:
 
 \newtheorem*{thm1}{Theorem 1}
\begin{thm1}
There exists a global $C^1$ solution of PDE (\ref{PDE}) on $[0,T]\times\mathbb{R}$. 
\end{thm1}

Finally, we establish uniqueness of solutions. Suppose $u_1(t,\theta)$ and $u_2(t,\theta)$ are $C^1$, solutions of (\ref{PDE}) and (\ref{simplemodelIC}) on $[0,T]\times\mathbb{R}$. Consider $\int_0^{2\pi}{\left((u_1(t,\theta)-u_2(t,\theta)\right)^2\ d\theta}$. Differentiating with respect to $t$ we find
{\fontsize{7}{6}
\begin{align}
    \frac{d}{dt}\int_{0}^{2\pi}&(u_1-u_2)^2(t,\theta)\ d\theta=2\int_{0}^{2\pi}(u_1-u_2)(t,\theta)\frac{\partial(u_1-u_2)}{\partial t}(t,\theta)\ d\theta\nonumber\\
    =\label{using_PDE}
    &-2k\int_{0}^{2\pi}(u_1-u_2)(t,\theta)\frac{\partial (u_1-u_2)}{\partial\theta}(t,\theta)\int_{0}^{2\pi}u_2(t,\hat{\theta})\sin(\hat{\theta}-\theta)d\hat{\theta}\ d\theta
    \nonumber\\
    &+2k\int_{0}^{2\pi}(u_1-u_2)(t,\theta)\frac{\partial u_1}{\partial \theta}(t,\theta)\int_0^{2\pi}(u_2-u_1)(t,\hat{\theta})\sin(\hat{\theta}-\theta)\ d\hat{\theta}\ d\theta
    \nonumber\\
    &+2k\int_{0}^{2\pi}(u_1-u_2)^2(t,\theta)\int_{0}^{2\pi}u_2(t,\hat{\theta})\cos(\hat{\theta}-\theta)\ d\hat{\theta}\ d\theta\nonumber\\
    &+2k\int_{0}^{2\pi}(u_1-u_2)(t,\theta)u_1(t,\theta)\int_{0}^{2\pi}(u_1-u_2)(t,\hat{\theta})\cos(\hat{\theta}-\theta)\ d\hat{\theta}\ d\theta\\
    \label{integral_by_parts_fin}=&\ 2k\int_{0}^{2\pi}(u_1-u_2)(t,\theta)\frac{\partial u_1}{\partial \theta}(t,\theta)\int_0^{2\pi}(u_2-u_1)(t,\hat{\theta})\sin(\hat{\theta}-\theta)\ d\hat{\theta}\ d\theta
    \nonumber\\
    &+k\int_{0}^{2\pi}(u_1-u_2)^2(t,\theta)\int_{0}^{2\pi}u_2(t,\hat{\theta})\cos(\hat{\theta}-\theta)\ d\hat{\theta}\ d\theta\nonumber\\
    &+2k\int_{0}^{2\pi}(u_1-u_2)(t,\theta)u_1(t,\theta)\int_{0}^{2\pi}(u_1-u_2)(t,\hat{\theta})\cos(\hat{\theta}-\theta)\ d\hat{\theta}\ d\theta\nonumber\\
    \leq& 2k\left(\sqrt{2\pi}\|\frac{\partial u_1}{\partial\theta}\|_{\infty}+ \pi\|u_2\|_{\infty}+\sqrt{2\pi}\|u_1\|_{\infty}\right)\int_{0}^{2\pi}(u_1-u_2)^2(t,\theta)\ d\theta, 
\end{align}
}
where $\|\cdot\|_{\infty}$ denotes the supremum on $[0,T]\times[0,2\pi]$, and where we obtain (\ref{integral_by_parts_fin}) after using integration by parts to combine the first and third integral in (\ref{using_PDE}). Hence, by Gronwall's inequality,
\begin{equation}
\int_{0}^{2\pi}(u_1-u_2)^2(t,\theta)\ d\theta \equiv \int_{0}^{2\pi}(u_1-u_2)^2(0,\theta)\ d\theta =0.
\end{equation}
We have established the following theorem. 
\newtheorem*{thm2}{Theorem 2}
\begin{thm2}
$C^1$ solutions of (\ref{PDE}) and (\ref{simplemodelIC}) on $[0,T]\times\mathbb{R}$ are unique. 
\end{thm2}

\section{Characterization of the solution}

In characterizing the solution density it is useful to note that the subsequence from the proof of Theorem 1 can also be chosen so that the characteristic curves, $y_n(t;\theta_0),$ converge uniformly. Indeed, previously we have seen that for $(t,\theta_0)\in [0,T]\times\mathbb{R}$, $|\dot{y}_n(t;\theta_0)|=|g_n(t;\theta_0)|\leq k$, and $\lvert\frac{\partial y_n}{\partial \theta_0}(t;\theta_0)\rvert\leq e^{kT}$. From these bounds, and since $y_n(0;\theta_0)=\theta_0$, we see that $\left\{y_n\right\}_{n=1}^{\infty}$ 
is uniformly bounded and equicontinuous on $[0,T]\times[0,2\pi].$ Moreover, since 
{\small\begin{align}
    \dot{y}_n(t,\theta_0)=g_n(t;\theta_0)=k\int_{0}^{2 \pi}\rho_{n-1}(t, \hat{\theta})\sin(\hat{\theta}-y_{n}(t;\theta_0))\ d\hat{\theta},\label{yndot}
\end{align}} 
and  
{\fontsize{8}{7}\begin{align}
\frac{\partial y_n}{\partial\theta_0}(t,\theta_0)=&\mathrm{exp}\left(-\int_0^t f_{n}(s;\theta_0)\ ds\right)\nonumber\\
=&\mathrm{exp}\left(-k\int_0^t{ \int_0^{2\pi}{\rho_{n-1}(s,\hat{\theta})\cos(\hat{\theta}-y_{n}(s;\theta_0))d\hat{\theta}}\ ds}\right),\label{pynptheta0}
\end{align}}
we see that $\left\{\dot{y}_n\right\}_{n=1}^{\infty}$ and $\left\{\frac{\partial y_n}{\partial\theta_0}\right\}_{n=1}^{\infty}$ are also equicontinuous. Indeed,
{\fontsize{8}{7}\begin{align}
 |\dot{y}_n(t_1;\theta_1)-\dot{y}_n(t_2;\theta_2)|=&|g_n(t_1;\theta_1)-g_n(t_2;\theta_2)|\label{uniform_cont_ydot}\nonumber\\
 =&\bigg\lvert k\int_{0}^{2 \pi}\rho_{n-1}(t_1, \hat{\theta})\sin(\hat{\theta}-y_{n}(t_1;\theta_1))\nonumber\\
 &-\rho_{n-1}(t_2, \hat{\theta})\sin(\hat{\theta}-y_{n}(t_2;\theta_2))\ d\hat{\theta}\bigg\rvert\nonumber\\
 \leq& \lvert y_{n}(t_1;\theta_1)-y_{n}(t_2;\theta_2)\rvert\nonumber \\
 &+\int_{0}^{2 \pi}{\lvert \rho_{n-1}(t_1, \hat{\theta})-\rho_{n-1}(t_2, \hat{\theta})\rvert\ d\hat{\theta}}
\end{align}}
Similarly, 
{\fontsize{8}{7}
\begin{align}
|f_n(t_1;\theta_1)-f_n(t_2;\theta_2)|
 \leq& \lvert y_{n}(t_1;\theta_1)-y_{n}(t_2;\theta_2))\rvert\nonumber\\
 &+ \int_{0}^{2 \pi}{\lvert \rho_{n-1}(t_1, \hat{\theta})-\rho_{n-1}(t_2, \hat{\theta})\rvert\ d\hat{\theta}}.
 \end{align}}
 Thus, the equicontinuity of $\left\{y_n\right\}_{n=1}^{\infty}$ and $\left\{\rho_n\right\}_{n=0}^{\infty}$ on $[0,T]\times[0,2\pi]$, implies that of $\left\{\dot{y}_n\right\}_{n=1}^{\infty}$ and $\left\{f_n\right\}_{n=1}^{\infty}.$
 Also, by the mean value theorem 
 {\fontsize{7}{6}\begin{align}
     \lvert\frac{\partial y_n}{\partial \theta_0}(t_1;\theta_1)-\frac{\partial y_n}{\partial \theta_0}(t_2;\theta_2)\rvert&\leq e^{KT}\left(\int_0^{t_1}{ |f_n(s,\theta_1)-f_n(s,\theta_2)| ds} + \int_{t_1}^{t_2}{ |f(s,\theta_2)| ds}\right)\label{uniform_cont_parypartheta0}\nonumber\\
     &\leq e^{KT}\left(\int_0^{T}{ |f_n(s,\theta_1)-f_n(s,\theta_2)| ds} + k|t_1-t_2|\right).
 \end{align}}
Hence, $\left\{\frac{\partial y_n}{\partial\theta_0}\right\}_{n=1}^{\infty}$ is also equicontinuous on $[0,T]\times\mathbb[0,2\pi]$. 
Therefore, we can choose a subsequence on which $y_n$, $\dot{y}_n$ and $\frac{\partial y_n}{\partial\theta_0}$ together with $\rho_n$ and its derivatives converge uniformly on $[0,T]\times[0,2\pi]$. We may also note that since $y_n(t,\theta_0+2\pi)=y_n(t,\theta_0)+2\pi,$ $y_n$ converges uniformly on $[0,T]\times\mathbb{R}.$ Passing the limit through (\ref{yndot}) and (\ref{pynptheta0}), we find:

\begin{equation}\label{ydot}
    \dot{y}(t;\theta_0)=k\int_{0}^{2 \pi}\rho(t, \hat{\theta})\sin(\hat{\theta}-y(t;\theta_0))\ d\hat{\theta}, 
    \end{equation} 
    
and
{\fontsize{8}{7}\begin{align}
\frac{\partial y}{\partial\theta_0}(t;\theta_0)=&\mathrm{exp}\left(-k\int_0^t{ \int_0^{2\pi}{\rho(s,\hat{\theta})\cos(\hat{\theta}-y(s;\theta_0))d\hat{\theta}}\ ds}\right)\label{pyptheta0}
\end{align}}

Moreover, there is a unique global solution $y(t;t_0,\theta)$ of \begin{align}
    \dot{y}(t)=k\int_{0}^{2 \pi}\rho(t, \hat{\theta})\sin(\hat{\theta}-y(t))\ d\hat{\theta},\end{align} subject to any initial data $y(t_0)=\theta$. By existence and uniqueness, if we define $\phi:\mathbb{R}^2\rightarrow\mathbb{R}^2$ by 
    \begin{equation}
    \phi(t,\theta_0):=(t,y(t;\theta_0))
    \end{equation}
    then $\phi$ is onto $\mathbb{R}^2$ and invertible with 
    \begin{equation}
    \phi^{-1}(t,\theta)=(t,y(0;t,\theta)).
    \end{equation}
    In fact, as in the proof of Corollary 1, we see that $\phi$ and $\phi^{-1}$ are $C^1$ with uniformly bounded derivatives on $[0,T]\times\mathbb{R}.$ Hence, $\phi$ and $\phi^{-1}$ are uniformly continuous on $[0,T]\times\mathbb{R}.$ 
    
    Now we will show that $\phi_n(t,\theta_0):=(t,y_{n+1}(t;\theta_0))$ and $\phi_n^{-1}(t,\theta):=(t,y_{n+1}(0;t,\theta))$ (see (\ref{geninitcondchar1-3_solution}) for clarification on this notation) converge uniformly to $\phi(t,\theta_0)$ and $\phi^{-1}(t,\theta)$, respectively, on $[0,T]\times\mathbb{R}.$
    \begin{proof}
    Uniform convergence of $y_n$ implies that of $\phi_n$. Since $\phi_n$ converges uniformly to $\phi$ on $[0,T]\times\mathbb{R}$, for every $\epsilon>0$, there exists $N\in\mathbb{N}$, so that for every $n>N$, and every $\theta\in\mathbb{R},$ 
    \begin{equation}
    \lvert\phi(\phi_n^{-1}(t,\theta))-(t,\theta)\rvert=\lvert\phi(\phi_n^{-1}(t,\theta))-\phi_n(\phi_n^{-1}(t,\theta))\rvert <\epsilon.
    \end{equation}
    That is, $\phi(\phi_n^{-1}(t,\theta))\rightarrow \phi(\phi^{-1}(t,\theta))$ uniformly as $n\rightarrow\infty$. Since $\phi^{-1}$ is uniformly continuous, given any $t\in[0,T]$, we can apply  $\phi^{-1}$ to this sequence to find that $\phi_n^{-1}(t,\theta)\rightarrow \phi^{-1}(t,\theta)$ uniformly on $[0,T]\times\mathbb{R}$.\end{proof}
    
    Now we show that uniform convergence of $\phi_n^{-1}$ together with uniform convergence and equicontinuity of $\left\{y_n(t;\theta_0)\right\}_{n=1}^{\infty}$ gives uniform convergence of $Y_n(s, t, \theta):=y_{n}(s;y_n(0;t,\theta))$ on $[0,T]\times [0,T]\times\mathbb{R}$.
    \begin{proof}
    Since $y_n$ is an equicontinuous sequence, given $\epsilon>0$, there exists an $\delta>0$, so that for every $s\in[0,T]$, $u_1, u_2\in\mathbb{R}$, and $m\in\mathbb{N}$, 
    \begin{equation}
    |u_1-u_2|<\delta,\rightarrow|y_{m}(s;u_1)-y_{m}(s;u_2)|<\frac{\epsilon}{2}.
    \end{equation}
    Since $y_n$  is a uniformly convergent sequence, there exists an $M>0,$ so that if $m,n>M$ then 
    \begin{equation}
    |y_{n}(s;u_1)-y_{m}(s;u_1)|<\frac{\epsilon}{2}.
    \end{equation}
    Also, since $\phi^{-1}_n$ is a uniformly convergent sequence, given $\delta>0$, there exists an $N$ so that for $n,m>N$, $t\in[0,T]$ and $\theta\in\mathbb{R},$
  \begin{equation}
  |y_n(0;t,\theta)-y_m(0;t,\theta)|<\delta.
  \end{equation}
  Thus, if $\delta$ and $\epsilon$ are as above, and $m,n>\max\left\{M,N\right\}$, then 
{\fontsize{7}{6}\begin{align}|y_{n}(s;y_n(0;t,\theta))-y_{m}(s;y_m(0;t,\theta))|&\leq |y_{n}(s;y_n(0;t,\theta))-y_{m}(s;y_n(0;t,\theta))|\\
&+|y_{m}(s;y_n(0;t,\theta))-y_{m}(s;y_m(0;t,\theta))|\\
  &<\epsilon.
  \end{align}}
  That is, $Y_n(s, t, \theta)$ is uniformly Cauchy, and hence uniformly convergent.\end{proof}
By the previous result, we may pull the limit through the solution formula (\ref{rhoSolution}), to find
{\fontsize{7}{6}\begin{eqnarray}
    \rho(t,\theta)=\rho_{0}(y(0;t,\theta)) \mathrm{exp}\left(k\int_{0}^{t}\int_{0}^{2\pi}\rho(s,\hat{\theta})\cos(\hat{\theta}-y(s;y(0;t,\theta)))\ d\hat{\theta}\ ds\right).
\end{eqnarray}}

This gives a characterization of the solution.
\newtheorem*{thm3}{Theorem 3}
\begin{thm3}
Continuously differentiable solutions of (\ref{simplifiedmodel})-(\ref{periodicity}) on $[0,T]\times\mathbb{R}$ are characterized by the following system:
{\fontsize{7}{6}\begin{align}
    \rho(t,\theta)&=\rho_{0}(y(0;t,\theta)) \mathrm{exp}\left(k\int_{0}^{t}\int_{0}^{2\pi}\rho(s,\hat{\theta})\cos(\hat{\theta}-y(s;y(0;t,\theta)))\ d\hat{\theta}\ ds\right),\\
    \dot{y}(t;\theta_0)&=k\int_{0}^{2 \pi}\rho(t, \hat{\theta})\sin(\hat{\theta}-y(t;\theta_0))\ d\hat{\theta}.\label{solution1}
\end{align}}
\end{thm3}
Note that since $T>0$ is arbitrary, the solution characterization is valid for all $T.$

\section{Some notes about the stable distribution.} 

Let 
\begin{equation}
g(\theta;t)=\int_0^{2\pi}\rho(\hat{\theta};t)\sin(\hat{\theta}-\theta)d\hat{\theta},\
\end{equation}
and 
\begin{equation}
f(\theta;t)=\int_0^{2\pi}\rho(\hat{\theta};t)\cos(\hat{\theta}-\theta)d\hat{\theta},
\end{equation}
where we treat $t$ as a fixed parameter.
Then 
\begin{equation}
g'(\theta;t)=-f(\theta;t)
\end{equation}
and 
\begin{equation}
g''(\theta;t)=-g(\theta;t).
\end{equation}
It follows that $g(\theta;t)=A(t)\cos(\theta)+B(t)\sin(\theta).$ Furthermore, if  $g(\theta;t)$ is not equivalent to $0$, there must exist exist two angles $\psi(t)\in[0,2\pi]$ so that $g(\psi(t);t)=0.$ Also, since $f^2(\theta;t)+g^2(\theta;t)\equiv C^2(t),$ we have, $C^2(t) = f^2(\psi(t);t).$ Finally, since $f(\psi(t)+\pi;t)=-f(\psi(t);t)$ we can choose $\psi(t)$ so that $f(\psi(t);t)=C(t)>0$. Thus, putting $x(\theta)=g(\psi(t)-\theta;t)$ and $y(\theta)=f(\psi(t)-\theta;t)$ we find that $x''(\theta)=-x(\theta),$ $x'(0)=f(\psi(t);t)=C(t),$ and $x(0)=g(\psi(t);t)=0.$ Thus, $x(\theta)=C(t)\sin(\theta)$ and $y(\theta)=C(t)\cos(\theta).$ That is, 
\begin{equation}
g(\theta;t)=C(t)\sin(\psi(t)-\theta)
\end{equation}and 
\begin{equation}
f(\theta;t)=C(t)\cos(\psi(t)-\theta).
\end{equation}

Now the differential equation for $y$ can be written as 
\begin{equation}
\dot{y}(t;\theta_0)=C(t)\sin(\psi(t)-y(t;\theta_0)).
\end{equation}
We would like to show that $y(t;\theta_0)\to\psi(t),$ $mod\ 2\pi,$ as $t\to\infty$.

In order to analyze the behavior of $y(t;\theta_0)$, we will first characterize the behavior of $C(t)$, which provides a measure of the order of the system. Indeed, $|C(t)|\leq 1,$ and $C(t)=1$ when $\rho(t,\theta)=\delta (\theta-\psi(t))$, where $\delta$ denotes the Dirac delta distribution. Note 
\begin{align}
C(t)=&\int_0^{2\pi}{\rho(t,\theta)\cos(\theta-\psi(t))d\theta}\\
=&\int_0^{2\pi}\rho_{0}(y(0;t,\theta)) e^{\left(k\int_{0}^{t}\int_{0}^{2\pi}\rho(s,\hat{\theta})\cos(\hat{\theta}-y(s;y(0;t,\theta)))\ d\hat{\theta}\ ds\right)}\nonumber\\
&\times\cos(\theta-\psi(t))\ d\theta
\end{align}
Letting $u=y(0;t,\theta),$ equivalently $\theta=y(t;u),$ we have
{\fontsize{8}{7}
\begin{equation}
\frac{d\theta}{du}=\frac{\partial y}{\partial u}(t;u)= \mathrm{exp}\left(-k\int_0^t{\int_0^{2\pi}{ \rho_{0}(t,\hat{\theta})\cos(\hat{\theta}-y(s;u))d\hat{\theta}}ds}\right).
\end{equation}} Hence, under this change of variable 
\begin{equation}
C(t)=\int_0^{2\pi}{\rho_{0}(u) \cos(\psi(t)-y(t;u))\ du}.
\end{equation}
Note here, we use the fact that $y(t;\theta_0+2\pi)=y(t;\theta_0)+2\pi$ and the fact that $\rho$ and $\cos$ are $2\pi$-periodic to preserve the domain of integration under the change of variable. Indeed, if $u_0=y(0;t,0)$ then $y(t;u_0)=0$, so $y(t;u_0+2\pi)=2\pi$, and $u_0+2\pi=y(0;t,2\pi)$. Choosing $k\in\mathbb{Z}$ so that $u_0\leq2k\pi\leq u_0+2\pi,$
{\begin{align}
C(t)=&\int_{u_0}^{u_0+2\pi}{\rho_{0}(u) \cos(\psi(t)-y(t;u))\ du}\\
=&\int_{u_0}^{2k\pi}{\rho_{0}(u) \cos(\psi(t)-y(t;u))\ du}\nonumber\\
&+\int_{2k\pi}^{u_0+2\pi}{\rho_{0}(u) \cos(\psi(t)-y(t;u))\ du}\\
=&\int_{u_0+2\pi}^{2k\pi+2\pi}{\rho_{0}(u-2\pi) \cos(\psi(t)-y(t;u-2\pi))\ du}\nonumber\\
&+\int_{2k\pi}^{u_0+2\pi}{\rho_{0}(u) \cos(\psi(t)-y(t;u))\ du}\\
=&\int_{u_0+2\pi}^{2k\pi+2\pi}{\rho_{0}(u) \cos(\psi(t)-y(t;u)+2\pi)\ du}\nonumber\\
&+\int_{2k\pi}^{u_0+2\pi}{\rho_{0}(u) \cos(\psi(t)-y(t;u))\ du}\\
=&\int_{u_0+2\pi}^{2k\pi+2\pi}{\rho_{0}(u) \cos(\psi(t)-y(t;u))\ du}\nonumber\\
&+\int_{2k\pi}^{u_0+2\pi}{\rho_{0}(u) \cos(\psi(t)-y(t;u))\ du}\\
=&\int_{2k\pi}^{2k\pi+2\pi}{\rho_{0}(u) \cos(\psi(t)-y(t;u))\ du}\\
=&\int_{0}^{2\pi}{\rho_{0}(u) \cos(\psi(t)-y(t;u))\ du}.
\label{change_of_domain}
\end{align}}
Under the same change of variable
{\fontsize{8}{7}
\begin{equation}
\int_0^{2\pi}\rho_0(u)\sin(y(t;u)-\psi(t))\ du = \int_0^{2\pi}\rho(t,\theta)\sin(\theta-\psi(t))\ d\theta = 0.
\end{equation}}
Hence,
{\fontsize{8}{7}\begin{align}
C'(t)&=-\int_0^{2\pi}\rho_0(u)\sin(\psi(t)-y(t;u))(\psi'(t)-C(t)k\sin(\psi(t)-y(t;u))\ du\nonumber\\
&=C(t)\int_0^{2\pi}k\rho_0(u)\sin^2(\psi(t)-y(t;u))\ du,
\end{align}}
and, $C'(t)\geq 0.$ Since also, $C(t)\leq 1,$ (recall $\rho_0$ is a probability density), $C(t)$ converges to a finite value, $C^*,$ as $t\rightarrow\infty$. Also, 
we find that 
{\fontsize{8}{7}\begin{equation}
   C(t)=C(0)\exp{\left(k\int_0^t{\int_0^{2\pi}{\rho_0(u)\sin^2(\psi(s)-y(s;u))\ du}\ ds}\right) }.\label{C}
\end{equation}}
So 
{\small\begin{equation}
\lim_{t\rightarrow\infty}k\int_0^t\int_0^{2\pi}\rho_0(u)\sin^2(\psi(s)-y(s;u))\ du\ ds=\ln\left(\frac{C^*}{C(0)}\right)
\label{finite_limit}\end{equation}}
is finite. Since, in addition, $k\int_0^{2\pi}\rho_0(u)\sin^2(\psi(t)-y(t;u))\ du$ is nonnegative and uniformly continuous on $[0,\infty)$ (Note its derivative is bounded.), 
\begin{equation}
\lim_{t\rightarrow\infty}k\int_0^{2\pi}\rho_0(u)\sin^2(\psi(t)-y(t;u))\; du=0.\end{equation}

Implicitly differentiating \begin{equation}
    \int_0^{2\pi}{\rho_0(u)\sin(y(t;u)-\psi(t))du}=0,
\end{equation} we find that $\psi(t)$ is differentiable, and 
{\fontsize{8}{7}\begin{equation}
\psi'(t) = k\int_0^{2\pi}\rho_0(u)\cos(\psi(t) - y(t; u))\sin(\psi(t)-y(t;u))\ du.
\end{equation}}
Hence,
{\fontsize{8}{7}
\begin{eqnarray}
\psi'(t)-y'(t;\theta_0) =&& -kC(t)\sin(\psi(t) - y(t;\theta_0)) \\
&&+ k\int_0^{2\pi}\rho_0(u)\cos(\psi(t) - y(t; u))\sin(\psi(t)-y(t;u))du\nonumber
\end{eqnarray}}
Define 
{\fontsize{9}{8}
\begin{equation}
\epsilon(t) := k\int_0^{2\pi}\rho_0(u)\cos(\psi(t) - y(t; u))\sin(\psi(t)-y(t;u))\ du.
\end{equation}}
Then,
\begin{eqnarray}
|\epsilon(t)| \leq&& k\left(\int_0^{2\pi}\rho_0(u)\cos^2(\psi(t) - y(t; u))\ du\right)^\frac{1}{2}\nonumber\\
&&\times\left(\int_0^{2\pi}\rho_0(u)\sin^2(\psi(t)-y(t;u)\ du\right)^\frac{1}{2}.
\end{eqnarray}
Since
\begin{equation}
\left(\int_0^{2\pi}\rho_0(u)\cos^2(\psi(t) - y(t; u))du\right)^\frac{1}{2} \leq 1
\end{equation}
and 
\begin{equation}
\lim_{t\rightarrow\infty}\left(\int_0^{2\pi}\rho_0(u)\sin^2(\psi(t) - y(t; u))du\right)^\frac{1}{2} =0
\end{equation}
\begin{equation}
\lim_{t\to\infty}\epsilon(t)=0.
\end{equation}
Setting 
\begin{equation}
\epsilon^*(T):=\sup\left\{|\epsilon(t)|:\   t\geq T\right\},
\end{equation}
we see that for $t>T$,
{\fontsize{9}{8}\begin{eqnarray}
 \psi'(t)- y'(t;\theta_0)&>&-kC(t)\sin(\psi(t)-y(t;\theta_0)) - \epsilon^*(T)  \\
\psi'(t)- y'(t;\theta_0)&<&-kC(t)\sin(\psi(t)-y(t;\theta_0)) + \epsilon^*(T).
\end{eqnarray}}
Since, in addition, $C(t)$ is monotone increasing and positive, 
\begin{itemize}
\item[(i)]$\psi(t) - y(t;\theta_0)$ is a subsolution of  
\begin{equation}
    f_+'(t) = -kC(T)\sin(f_+(t)) + \epsilon^*(T),\label{upperODE}
\end{equation}
for $t>T$, and $0< \psi(t) - y(t;\theta_0) < \pi$.\vspace{.5 cm}
\item[(ii)]$\psi(t) - y(t;\theta_0)$ is a supersolution of 
\begin{equation}
    f'_-(t) = -kC(T) \sin(f_-(t))-\epsilon^*(T),\label{lowerODE}
\end{equation}
for $t>T$ and $\pi< \psi(t) - y(t;\theta_0) < 2\pi$.\vspace{.5 cm}
\end{itemize}
(The previous claim involves a slight extension of Theorem 1.2 of \cite{teschl}. Please see the appendix.)

Note that since $\epsilon^*(t)\downarrow0$, and $C(t)$ is increasing, for $T$ large (\ref{upperODE}) has stable and unstable steady-state solutions, $s^T_++2j\pi$ and $u^T_++2j\pi,$ $j\in\mathbb{Z}$, defined by 
\begin{equation}
\sin(s_+^T) = \frac{\epsilon^*(T)}{kC(T)},\quad 0\leq s_+^T \leq \frac{\pi}{2}
\end{equation}
and 
\begin{equation}
\sin(u_+^T) = \frac{\epsilon^*(T)}{kC(T)},\quad \frac{\pi}{2}\leq u_+^T \leq \pi.
\end{equation}
Similarly, (\ref{lowerODE}) has stable and unstable steady-state solutions, $u_-^T+2j\pi$ and $s_-^T+2j\pi$, defined by 
\begin{equation}
\sin(u^T_-) = -\frac{\epsilon^*(T)}{kC(T)},\quad \pi\leq u_-^T \leq \frac{3\pi}{2},
\end{equation}
and 
\begin{equation}
\sin(s_-^T) = -\frac{\epsilon^*(T)}{kC(T)},\quad \frac{3\pi}{2}\leq s_-^T \leq 2\pi.
\end{equation}
It follows that if 
\begin{equation}
s_-^T+2(j-1)\pi<\psi(t) - y(t;\theta_0)<s_+^T+2j\pi
\end{equation}
for some $t>T,$ then 
\begin{equation}
s_-^T+2(j-1)\pi<\psi(t) - y(t;\theta_0)<s_+^T+2j\pi
\end{equation}
for all later times. 

Also, if for some $t>T,$ 
\begin{equation}
s_+^T+2j\pi<\psi(t) - y(t;\theta_0)<u_+^T+2j\pi,
\end{equation}
then 
\begin{equation}
s_-^T+2(j-1)\pi\leq \limsup\left\{\psi(t) - y(t;\theta_0)\right\} \leq s_+^T+2j\pi,
\end{equation}
and 
\begin{equation}
s_-^T+2(j-1)\pi\leq \liminf\left\{\psi(t) - y(t;\theta_0)\right\} \leq s_+^T+2j\pi.
\end{equation}
Similarly, if 
\begin{equation}
u_-^T+2j\pi<\psi(t) - y(t;\theta_0)<s_-^T+2j\pi,
\end{equation}
then 
\begin{equation}
s_-^T+2j\pi\leq \limsup\left\{\psi(t) - y(t;\theta_0)\right\} \leq s_+^T+2(j+1)\pi,
\end{equation}
and 
\begin{equation}
s_-^T+2j\pi\leq \liminf\left\{\psi(t) - y(t;\theta_0)\right\} \leq s_+^T+2(j+1)\pi.
\end{equation}

That is, if there exists $t>T$ so that $\psi(t) - y(t;\theta_0)\notin\bigcup_{i\in\mathbb{N}}[u_+^{T}+2i\pi, u_-^{T}+2i\pi]$, then there exists $j\in\mathbb{N}$ so that $\psi(t) - y(t;\theta_0)$ is \textit{asymptotically contained} in $[s_-^{T}+2(j-1)\pi, s_+^{T}+2j\pi]$. This means that if $U$ is an open set with $[s_-^{T}+2(j-1)\pi, s_+^{T}+2j\pi]\subset U$, then there exists $\tau$, so that for $t>\tau$, $\psi(t) - y(t;\theta_0)\in U$. 

Furthermore, since $u_+^{T}$ is monotone increasing in $T,$ and $u_-^{T}$ is monotone decreasing in $T$, we see that if $\psi(t) - y(t;\theta_0)$ is asymptotically contained in $[s_-^{T}+2(j-1)\pi, s_+^{T}+2j\pi]$, then $\psi(t) - y(t;\theta_0)$ is asymptotically contained in $[s_-^{\tau}+2(j-1)\pi, s_+^{\tau}+2j\pi]$, for $\tau>T.$ Indeed, 
{\fontsize{9}{8}\begin{align}
[s_-^{T}+2(j-1)\pi, s_+^{T}+2j\pi]
&\subset\left(\bigcup_{i\in\mathbb{N}}[u_+^{T}+2i\pi, u_-^{T}+2i\pi]\right)^C\nonumber\\
&\subset\left(\bigcup_{i\in\mathbb{N}}[u_+^{\tau}+2i\pi, u_-^{\tau}+2i\pi]\right)^C.
\end{align}}
And, since $\left(\bigcup_{i\in\mathbb{N}}[u_+^{\tau}+2i\pi, u_-^{\tau}+2i\pi]\right)^C$ is open, if $\psi(t) - y(t;\theta_0)$ is asymptotically contained in $[s_-^{T}+2(j-1)\pi, s_+^{T}+2j\pi],$ then there exists $t>\tau$, so that $\psi(t) - y(t;\theta_0)\notin\bigcup_{i\in\mathbb{N}}[u_+^{\tau}+2i\pi, u_-^{\tau}+2i\pi]$. Hence, $\psi(t) - y(t;\theta_0)$ is asymptotically contained in $[s_-^{\tau}+2(j-1)\pi, s_+^{\tau}+j\pi].$ 

Finally, since $s_{-}^T+2(j-1)\pi\rightarrow 2j\pi$ and $s_{+}^T+2(j)\pi\rightarrow 2j\pi$, as $T\to\infty$, we arrive at the following lemma:

\newtheorem*{lem8}{Lemma 8}
\begin{lem8}
If there exists $T$ and $t>T$ so that $\psi(t) - y(t;\theta_0)\notin\bigcup_{i\in\mathbb{N}}[u_+^{T}+2i\pi, u_-^{T}+2i\pi],$ then there exists $j\in\mathbb{Z}$ so that $\psi(t) - y(t;\theta_0)\rightarrow 2j\pi$ as $t\rightarrow\infty$.
\end{lem8}

\newtheorem*{lem9}{Lemma 9}
\begin{lem9}
Fix $T>0$. If for all $t > T$, $\psi(t) - y(t;\theta_0) \in \bigcup_{j\in\mathbb{Z}}[u_+^T+2j\pi, u_-^T+2j\pi],$ then there exists $j\in\mathbb{Z}$ so that $\psi(t) - y(t;\theta_0) \to (2j+1)\pi$ as $t\to\infty.$
\end{lem9}

\begin{proof} By continuity, if $\psi(t) - y(t;\theta_0) \in \bigcup_{k\in\mathbb{Z}}[u_+^T+2j\pi, u_-^T+2j\pi]$ for all $t>T$, then there exists $j\in\mathbb{Z}$ so that $\psi(t) - y(t;\theta_0) \in [u_+^T+2j\pi, u_-^T+2j\pi]$ for all $t>T$. Without loss of generality, let $j=0$, so that $\psi(t) - y(t;\theta_0) \in [u_+^{T}, u_-^{T}]$ for all $t > T$. Suppose toward a contradiction that $\lim\limits_{t\to\infty}\psi(t) - y(t;\theta_0)) \neq \pi$. Then there exists $\epsilon>0$, so that for every $\tau$, there exists $t^*>\tau$ so that 
\begin{equation}
|\psi(t) - y(t^*;\theta_0)-\pi| > \epsilon.
\end{equation}
Since $u_-^{\tau}$ and $u_+^{\tau}$ converge to $\pi$ as $\tau$ approaches infinity, there exists $\tau>T$ so that 
\begin{equation}
|u_-^{\tau}-\pi|=|u_+^{\tau}-\pi|<\frac{\epsilon}{2}.
\end{equation}
Hence, there exists $t^*>\tau>T$ so that
\begin{equation}\label{notin_ut}\psi(t^*) - y(t^*;\theta_0)\notin [u_+^{\tau}, u_-^{\tau}].\end{equation} Since $\psi(t) - y(t;\theta_0)\in [u_+^{T}, u_-^{T}]$ for all $t>T$, this means that 
\begin{equation}\label{notin_union_ut}\psi(t^*) - y(t^*;\theta_0)\notin\bigcup_{i\in\mathbb{N}}[u_+^{\tau}+2i\pi, u_-^{\tau}+2i\pi].\end{equation}
Hence $\psi(t) - y(t;\theta_0)$ is asymptotically contained in $[s_-^{\tau}, s_+^{\tau}]\subset[s_-^T, s_+^T]\subset[u_+^T, u_-^T]^C$, contradicting that $\psi(t) - y(t;\theta_0) \in [u_+^T, u_-^T]$ for all $t > T.$ 
\end{proof}

\newtheorem*{lem10}{Lemma 10}
\begin{lem10}
There exists at most one initial value $\theta_0\in[0, 2\pi)$ so that $\psi(t) - y(t;\theta_0)\to\pi,\mod 2\pi,$ as $t\to\infty.$
\end{lem10}

\begin{proof}
Suppose there exits $\theta_1, \theta_2\in[0,2\pi)$ with $\theta_1<\theta_2$, and there exists $i,j\in\mathbb{Z}$ so that $\lim\limits_{t\to\infty} y(t;\theta_1)-\psi(t) = (2j+1)\pi$ and $\lim\limits_{t\to\infty} y(t;\theta_2)-\psi(t)=(2i+1)\pi.$ Since  $0\leq\theta_1<\theta_2<2\pi$ and $y(t;\theta+2\pi)=y(t;\theta)+2\pi$, by uniqueness of solutions, for all $t$, 
\begin{equation}
y(t;0)\leq y(t;\theta_1)<y(t;\theta_2)<y(t;0)+2\pi.
\end{equation}
Hence $i=j$, and there exists $T$, so that for $t>T$, 
{\fontsize{8}{7}
\begin{equation}
(2j+1)\pi-\frac{\pi}{4} < y(t;\theta_1)-\psi(t) < y(t;\theta_2)-\psi(t) < (2j+1)\pi+\frac{\pi}{4}.
\end{equation}}
Thus, for $t>T$,
{\fontsize{8}{7}\begin{eqnarray}
(y(t,\theta_2)-y(t,\theta_1))' &=&-kC(t)(\sin(y(t,\theta_2)-\psi(t))-\sin(y(t,\theta_1)-\psi(t))\nonumber\\
&=& -kC(t)\cos(\theta(t))(y(t,\theta_2)-y(t,\theta_1));\\
&&\text{where } \theta(t)\in(y(t;\theta_1)-\psi(t),y(t,\theta_2)-\psi(t))\nonumber\\
&\geq&kC(t)\frac{\sqrt2}{2}(y(t;\theta_2)-y(t,\theta_1))\\
&\geq&kC(T)\frac{\sqrt2}{2}(y(t,\theta_2)-y(t,\theta_1))
\end{eqnarray}}
Putting $R=kC(T)\frac{\sqrt2}{2}$, we see that for $t>T$, $y(t,\theta_2) - y(t,\theta_1) > (y(T,\theta_2) -y(T,\theta_1))e^{R(t-T)}$, contradicting that $\lim\limits_{t\to\infty} y(t,\theta_1) = \lim\limits_{t\to\infty} y(t,\theta_2).$
\end{proof}

\newtheorem*{lem11}{Lemma 11}
\begin{lem11}
There exists a $\theta_0^C\in[0, 2\pi]$, and $j_c\in\mathbb{N},$ so that if $\theta_0\in[0, 2\pi]$ and $\theta_0<\theta_0^C,$ then $\psi(t) - y(t;\theta_0)\to 2j_C\pi$, and if $\theta_0\in[0, 2\pi]$ and $\theta_0>\theta_0^C$ then $\psi(t) - y(t;\theta_0)\to 2(j_C+1)\pi.$
\label{Lemma11}
\end{lem11}
\begin{proof}
From the contrapositive of Lemma 9 and Lemma 8, for each $\theta_0\in\mathbb{R}$, $y(t;\theta_0)-\psi(t)$ converges to $\pi, \mod 2\pi$ or to $0, \mod 2\pi.$ Note that continuity together with convergence to a limiting angle, mod $2\pi$, gives convergence to a limiting angle. Also, since $y(t;\theta_0+2\pi)=y(t;\theta_0)+2\pi,$ by uniqueness of solutions, there exists $j_c\in\mathbb{N},$ so that for each $\theta_0\in[0, 2\pi],$ $\lim_{t\to\infty}y(t;\theta_0)-\psi(t)\in[2j_c, 2(j_c+1)\pi].$
Let 
{\fontsize{8}{7}
\begin{equation}
\theta_+:=\inf\left\{\theta_0 :\ \theta_0\in[0,2\pi], \lim_{t\to\infty}y(t;\theta_0)-\psi(t)= 2(j_C+1)\pi\right\},
\end{equation}}
and
{\fontsize{8}{7}
\begin{equation}
\theta_-:=\sup\left\{\theta_0 :\ \theta_0\in[0,2\pi], \lim_{t\to\infty}y(t;\theta_0)-\psi(t)= 2j_C\pi\right\}.
\end{equation}}
By uniqueness, $\theta_-\leq\theta_+.$
If $\theta_-<\theta_+,$ then there exits $\theta_1, \theta_2$ with $\theta_-<\theta_1<\theta_2<\theta_+.$ Hence, $y(t;\theta_1)$ and $y(t;\theta_2)$ both converge to $(2j_C+1)\pi$, contradicting Lemma 10. Therefore, it must be that $\theta_-=\theta_+,$ so $\theta_0^C=\theta_-=\theta_+$ is the desired quantity.
\end{proof}

Since $y(t;\theta+2\pi)=y(t;\theta)+2\pi$, we have the following lemma.

\newtheorem*{lem12}{Lemma 12}
\begin{lem12}
There exists a unique $\theta_0^C\in[0, 2\pi)$, so that if $\theta_0\neq\theta_0^C,\mod{2\pi}$, then $\psi(t) - y(t;\theta_0)\to 0, \mod 2\pi$ as $t\to\infty.$
\label{Lemma12}
\end{lem12}

Now we can characterize the asymptotic behavior of $C(t).$

\newtheorem*{lem13}{Lemma 13}
\begin{lem13}
$C(t) \to 1$ as $t\to\infty$.
\label{Lemma13}
\end{lem13}

\begin{proof}
Let $\mu$ denote the measure induced by $\rho_0$ so that given $S\subseteq\mathbb{R}$, 
\begin{equation}
\mu(S)=\int_{S}{\rho_0(\theta)d\theta}.
\end{equation}
Let $j_C\in\mathbb{Z}$ be such that $y(t;\theta_0^C)-\psi(t)\to(2j_C+1)\pi$ as $t\to\infty.$  

By continuity of $\rho_0$, given $1>>\epsilon>0,$ we can choose $\delta$ so that $1>>\delta>0,$ and $\mu(B(\theta_0^C,\delta))<\frac{\epsilon}{4}.$ Note that since 
$\delta<2\pi$, 
\begin{equation}
y(t;\theta_0^C +\delta)-\psi(t)\to 2(j_C+1)\pi,
\end{equation}
and 
\begin{equation}
y(t;\theta_0^C -\delta)-\psi(t)\to 2j_C\pi.
\end{equation}
Thus, we can choose $T$ so that for $t \geq T$, 
{\fontsize{9}{8}\begin{equation}
2(j_C+1)\pi-\frac{\epsilon}{4}<y(t,\theta_0^C +\delta)-\psi(t)< 2(j_C+1)\pi+\frac{\epsilon}{4}
\end{equation}}
and  
\begin{equation}
2j_C\pi-\frac{\epsilon}{4}<y(t;\theta_0^C -\delta)-\psi(t)< 2j_C\pi+\frac{\epsilon}{4}. 
\end{equation}
If $\theta_0\in[\theta_0^C-\pi,\theta_0^C-\delta]\subset[\theta_0^C-2\pi+\delta, \theta_0^C-\delta],$ then 
{\fontsize{8}{7}\begin{equation}
y(t;\theta_0^C+\delta)-2\pi=y(t;\theta_0^C-2\pi+\delta)<y(t;\theta_0)<y(t;\theta_0^C-\delta),
\end{equation}}
and hence, for $t>T,$
\begin{equation}
2j_C\pi-\frac{\epsilon}{4}<y(t;\theta_0)-\psi(t)<2j_C\pi+\frac{\epsilon}{4}.
\end{equation}
That is, if $\theta_0\in[\theta_0^C-\pi,\theta_0^C-\delta],$ then for $t>T$, 
\begin{equation}
|y(t;\theta_0)-\psi(t)|<\frac{\epsilon}{4}, \mod{2\pi}.
\end{equation}

Similarly, if $\theta_0\in[\theta_0^C+\delta, \theta_0^C+\pi]\subset[\theta_0^C+\delta, \theta_0^C+2\pi-\delta],$ then for $t>T$, 
\begin{equation}
|y(t;\theta_0)-\psi(t)|<\frac{\epsilon}{4} \mod{2\pi}.
\end{equation}
It follows that if $\theta_0 \in B(\theta_0^C,\delta)$ and $t> T$, then 
\begin{eqnarray}
\label{e1}| y(t;\theta_0)-\psi(t) | &<& \frac{\epsilon}{4}\mod{2\pi},\\ \label{e2}|\sin(y(t;\theta_0)-\psi(t))| &<& \frac{\epsilon}{4},
\end{eqnarray} 
and 
\begin{equation}\label{e3}\cos( y(t;\theta_0)-\psi(t)) > 1-\frac{\epsilon}{4}.
\end{equation}
Finally, note that \begin{equation}\label{e4}\mu(\left(B(\theta_0^C,\delta)\right)^C\cap[0,2\pi])>1-\frac{\epsilon}{4}.\end{equation}
From the (\ref{e2})-(\ref{e4}), we see that for $t>T,$
\begin{align}
C(t) =& \int_{\theta_0^C-\delta}^{\theta_0^C+\delta} {\rho_0(u)\cos(\psi(t)-y(t;u))\ du}\nonumber\\
&+ \int_{\theta_0^C-\pi}^{\theta_0^C-\delta}{\rho_0(u)\cos(\psi(t)-y(t;u))\ du}\nonumber\\
&+ \int^{\theta_0^C+\pi}_{\theta_0^C+\delta}{\rho_0(u)\cos(\psi(t)-y(t;u))\ du}\\
\geq& -\frac{\epsilon}{4}+\left(1-\frac{\epsilon}{4}\right)^2\\
\geq& 1-\epsilon,
\end{align}
where we have used (\ref{change_of_domain}) to change the domain of integration in the expression for $C(t).$
Since $\epsilon>0$ was arbitrary, $C(t) \to 1$ as $t\to\infty.$
\end{proof}

Now we will show that $\psi(t)$ converges to a fixed value. First, note that since
{\fontsize{8}{7}
\begin{equation}
\frac{d}{dt}\int_0^{2\pi}\rho_0(u)y(t;u)\ du = kC(t)\int_0^{2\pi}\rho_0(u)\sin(\psi(t)-y(t;u))\ du\equiv 0,
\end{equation}}
{\fontsize{8}{7}
\begin{equation}\label{expect_y}
\int_0^{2\pi}\rho_0(u)y(t;u)\ du \equiv \int_0^{2\pi}\rho_0(u)y(0;u)\ du=\int_0^{2\pi}\rho_0(u)u\ du.
\end{equation}} 
Also,
{\fontsize{9}{8}
\begin{eqnarray}\label{lim_expect_y-psi}
\lim_{t\to\infty}\int_0^{2\pi}\rho_0(u)(y(t;u)-\psi(t))\ du=&&2j_C\pi\mu([0,\theta_0^C))\nonumber\\
&&+2(j_C+1)\pi\mu((\theta_0^C, 2\pi])\nonumber\\:=&&L.
\end{eqnarray}}
Indeed, since for each $u\in\mathbb{R}$, $y(t;u)-\psi(t)$ converges as $t\to\infty$, for each $u\in\mathbb{R}$, $\left\{y(t;u)-\psi(t):\ t>0\right\}$ is bounded. Moreover, since for every $u\in[0,2\pi]$, $y(t;u)\in[y(t;0), y(t;2\pi)]$, $\left\{y(t;u)-\psi(t):\ t>0, u\in[0,2\pi]\right\}$ is also bounded. Hence, there exists $M$ so that, for each $t>0$ and $u\in[0,2\pi],$ $|y(t;u)-\psi(t)|<M$.  Since $\rho_0$ is continuous, given $1>>\epsilon>0,$ we can choose $\delta$ so that $1>>\delta>0,$ and $\mu(B(\theta_0^C,\delta))<\frac{\epsilon}{3M}.$
Also, since for $u\in[0,\theta_0^C)$, $y(t;u)-\psi(t)\to 2j_C\pi$ as $t\to\infty$, and for $u\in(\theta_0^C,2\pi]$, $ y(t;u)-\psi(t)\to 2(j_C+1)\pi$ as $t\to\infty,$ there exists $T$, so that
\begin{itemize}
    \item[(i)] for $t>T,$ and $u\in[0,\theta_0^C-\delta]$, 
    \begin{equation}
    |y(t;u)-\psi(t)-2j_C\pi|<\frac{\epsilon}{3};
    \end{equation}
    \item[(ii)]for $t>T,$ and $u\in[\theta_0+\delta, 2\pi]$, 
    \begin{equation}
    |y(t;u)-\psi(t)-2(j_C+1)\pi|<\frac{\epsilon}{3}.
    \end{equation}
\end{itemize}
Hence for $t>T,$ and $\theta_0^C\in(0,2\pi)$
{\fontsize{8}{7}
\begin{eqnarray}\label{expect_y-psi_large_t}
\lvert\int_0^{2\pi}\rho_0(u)(y(t;u)&-&\psi(t))\ du - L\rvert<\nonumber\\
&&\int_0^{\theta_0^C-\delta}\rho_0(u)|y(t;u)-\psi(t)-2j_c\pi|\ du\nonumber\\
&+&\int_{\theta_0^C-\delta}^{\theta_0^C+\delta}\rho_0(u)|y(t;u)-\psi(t)|\ du\nonumber\\
&+&\int_{\theta_0^C+\delta}^{2\pi}\rho_0(u)|y(t;u)-\psi(t)-2(j_c+1)\pi|\ du\nonumber\\
&<&\epsilon.
\end{eqnarray}}
A similar equality holds for $\theta_0^C\in\left\{0,2\pi\right\}$. This establishes (\ref{lim_expect_y-psi}). 
From (\ref{expect_y}) and (\ref{lim_expect_y-psi}), we have the following lemma.
\newtheorem*{lem14}{Lemma 14}
\begin{lem14}
\begin{equation}
\lim_{t\to\infty}\psi(t)=\int_0^{2\pi}\rho_0(u)u\ du-L.
\end{equation}
\label{Lemma14}
\end{lem14}

\begin{proof}
{\fontsize{9}{8}
\begin{eqnarray}
     \psi(t)&=&\int_0^{2\pi}\rho_0(u)\psi(t)\ du\nonumber\\
     &=&-\int_0^{2\pi}\rho_0(u)(y(t;u)-\psi(t))\ du\nonumber + \int_0^{2\pi}\rho_0(u)y(t;u)\ du
\end{eqnarray}}
\end{proof}

The previous work also gives an alternate characterization of the solution. 

\newtheorem*{thm4}{Theorem 4}
\begin{thm4}
Continuously differentiable solutions of (\ref{simplifiedmodel})-(\ref{periodicity}) on $[0,T]\times\mathbb{R}$ are characterized by the following system:
{\fontsize{7}{6}\begin{align}
    \rho(t,\theta)&=\rho_{0}(y(0;t,\theta)) \mathrm{exp}\left(k\int_{0}^{t}\int_{0}^{2\pi}\rho_0(u)\cos(y(t;u)-y(s;y(0;t,\theta)))\ du\ ds\right),\label{solution2a}\\
    \dot{y}(t;\theta_0)&=k\int_{0}^{2 \pi}\rho_0(u)\sin(y(t;u)-y(t;\theta_0))\ du.\label{solution2b}
\end{align}}
\end{thm4}
Since $T>0$ is arbitrary, the solution characterization is valid for all positive $T.$
The utility of the characterization (\ref{solution2a})-(\ref{solution2b}), is that the projected characteristic equation can now be solved independently of the density equation. Hence, this characterization suggests a numerical algorithm for solving the model. Finally note that although the solution characterization was derived assuming a twice continuously differentiable initial density, the integral equations may be solvable under weaker constraints. That is (\ref{solution2a})-(\ref{solution2b}) can be investigated as a weak solution of the Kuramoto model with constant natural velocities. 

\section{Conclusion}
In this paper we used an iterative method combined with the method of characteristics for first order partial differential equations to show existence of global solutions to the continuum version of the Kuramoto model with identical oscillators. Our method yields a continuously differential, global solution and provides a characterization of the oscillator density in terms of a limiting projected characteristic, $y(t;\theta_0)$. This projected characteristic can be considered to describe the dynamics of oscillators with an initial phase of $\theta_0$. The solution characterization provides for an asymptotic analysis of the model solutions. We verify that, in the case of identical oscillators, for almost every $\theta_0$, $y(t;\theta_0)$ converges, mod $2\pi,$ to the average phase, $\psi(t)$. Moreover, the oscillator system asymptotically approaches perfect order, that is $C(t)\to1$ as $t\to\infty$. Our results open the door to considerable future work. The solution characterization suggests both a numerical algorithm for solving the model. In addition, we anticipate the methods developed here may be extendable to the continuum version of Kuramoto model with nonidentical oscillators.      

\flushleft{\textbf{Data Availability Statement:}
Data sharing is not applicable to this article as no new data were created or analyzed in this study.}

\appendix
\section{Appendix}
Here we provide a corollary of theorem 1.2 of \cite{teschl}, which is useful in section 4. The proof follows that in \cite{teschl}.

\newtheorem*{lem15}{Corollary A1}
\begin{lem15}
Let $y$ and $y_+$ satisfy, $y'(t)=F(t,y),$ $y_+'(t)=F_+(t,y_+),$ respectively. If there exists $t_0\in\mathbb{R}$ and $(a,b)\subseteq\mathbb{R} $ so that 
\begin{itemize}
    \item $y(t_0),y_+(t_0)\in(a,b)$,
    \item $y(t_0)\leq y_+(t_0),$ and
    \item $F(t,y)<F_+(t,y)$ for $y\in(a,b),$
\end{itemize}  
then $y(t)< y_+(t)$ while $y(t)\in(a,b).$ 
\label{Lemma15}
\end{lem15}

\begin{proof} 
Let $y, y_+, F, F_+, t_0,$ and $(a,b)$ satisfy the hypotheses of the corollary statement. Note that by Lemma 1.2 of \cite{teschl} if $y_+(t), y(t)\in(a,b)$ for $t_0\leq t\leq \hat{T}$, then $y_+$ is a supersolution of $y'(t)=F(t,y),$ on $[t_0,\hat{T})$, and hence, $y(t) < y_+(t)$ for $t\in(t_0,  \hat{T}$). Also, by continuity, there exists $\hat{T}$ so that $t_0<\hat{T},$ and $y(t), y_+(t) \in(a,b)$ for $t_0<t<\hat{T}.$ Hence, there exists $\hat{T}$ so that $y(t)<y_+(t)$ for $t_0<t<\hat{T}$

Suppose there exists $t_1 > t_0$ so that $y(t)\in(a,b)$ for $t_0\leq t \leq t_1,$ and yet $y(t_1)>y_+(t_1).$ Let \begin{equation}
E:=\left\{t:\ t_0<t\leq t_1, y(t)-y_+(t)\geq 0\right\},
\end{equation} and let $s=\inf(E).$
By the approximation property of infimum and the previous paragraph, $s>t_0$. Hence, $s\in(t_0, t_1]$, and, by continuity, $y_+(s)=y(s)\in(a,b)$. Then, by the third assumption in the corollary statement $\dot{y-y_+}(s)<0,$ contradicting that $s$ is the greatest lower bound of $E$.     
\end{proof}

\bibliography{main}

\end{document}